\documentclass[preprint,12pt]{elsarticle}

\usepackage{lineno,hyperref}
\usepackage{amstext}
\usepackage{graphicx}
\usepackage{epstopdf}
\usepackage{amssymb}
\usepackage{mathrsfs}
\usepackage{amsmath}
\usepackage{amsthm}
\usepackage{amsfonts}
\usepackage{framed}
\usepackage{verbatim}
\newcommand{\reals}{\mathbb{R}}
\newcommand{\DD}{\displaystyle}
\modulolinenumbers[5]

\journal{Journal of Computational Physics}

\begin{document}

\begin{frontmatter}

\title{Fast Algorithm for Simulating Lipid Vesicle Deformation I: Spherical Harmonic Approximation}


\author[mymainaddress]{Michael Mikucki}
\ead{mikucki@math.colostate.edu}

\author[mymainaddress]{Y. C. Zhou\corref{mycorrespondingauthor}}
\cortext[mycorrespondingauthor]{Corresponding author}
\ead{yzhou@math.colostate.edu}

\address[mymainaddress]{Department of Mathematics, Colorado State University, Fort Collins, CO 80523-1874}

\begin{abstract}
Lipid vesicles appear ubiquitously in biological systems. Understanding how the mechanical and 
intermolecular interations deform vesicle membrane is a fundamental question in biophysics. In this 
article we developed a fast algorithm to compute the surface configurations of lipid vesicles by 
introducing the surface harmonic functions to approximate the surfaces. This parameterization of the surfaces 
allows an analytical computation of the membrane curvature energy and its gradient for the efficient 
minimization of the curvature energy using a nonlinear conjugate gradient method. Our approach drastically 
reduces the degrees of freedom for approximating the membrane surfaces compared to the previously developed 
finite element and finite difference methods. Vesicle deformations with a reduced volume larger than
0.65 can be well approximated by using as small as 49 surface harmonic functions. The method thus has 
a great potential to reduce the computational expense of tracking multiple vesicles which deform for 
their interaction with external fields.
\end{abstract}


\end{frontmatter}


\section{Introduction}

This paper describes a fast numerical algorithm for computing the configuration
of lipid bilayer vesicles. Lipid bilayers are crucial components to
living systems. Being amphiphilic, lipid molecules have charged or polar 
hydrophilic head groups and hydrophobic tails. This allows lipids
in aqueous solution to aggregate into structures that entropically favor the
alignment of hydrophobic tails and the exposure of hydrophilic head groups
to water. A lipid bilayer is formed from the self-assembly of hydrophobic tails
of the two complementary layers. The bilayer will the then close so the hydrophobic core will not be exposed at the free edges, forming membranes
of cells and sub-celluar organelles. These membranes are semi-permeable boundaries separating the enclosure
from the surrounding environment. At the microscale, lipid bilayer membranes regulate the transportation of ions, proteins 
and other molecules between separated domains, and provide a flexible platform on which molecules can aggregate
to carry out vital chemical or physical reactions. 
At the mesoscale, for example, membranes of the red blood cells (RBCs) 
suspended in blood flow change their shape in response to the local flow conditions, and this change will in turn 
affect the RBC's ability of oxygen transport and the hydrodynamic properties of the blood flow \cite{EggletonC1998a,BagchiP2007a}. 
Deformation of the bilayer membranes can be driven by various types of force. At the microscale, driven forces are mainly the 
results of protein-membrane or membrane-membrane interactions, such as protein binding or insertion, lipid insertion or 
translation, and ubiquitous electrostatic interactions, to name a few \cite{Farsad2003}. At the mesoscale, hydrodynamic 
forces usually dominate \cite{WeiG2010a,KammR2002a}. Determination of the membrane geometry in response to 
protein-membrane, membrane-membrane, or fluid-membrane interactions is necessary for elucidating the structure-function 
relation of these interacted biological systems. 

The variation of lipid bilayer configurations can be characterized by its deformation energy.
This energy is the handle of almost all computational methods. Classical strain energy can be defined 
for lipid bilayers so their deformation can be described as elastic plates \cite{Zhou2010}. A crucial difference 
between the plates and bilayer membranes is missing in such models: a flat membrane can be subject to some shear 
deformation with zero energy cost provided the deformation is so slow that the viscous effect of lipids is negligible \cite{PowersT2005a}.
As a result, the deformation energy of a lipid bilayer is mostly attributed to the bending energy of the monolayers. 
The classical energy forms proposed by Canham \cite{Canham1970}, Helfrich \cite{Helfrich1973}, and Evans \cite{Evans1974}, are of this type, 
in which the deformation energy is defined be a quadratic function of the principle curvature of the surface. Other
components of the energy, such as those corresponding to the area expansion and contraction of the monolayers, 
and the osmotic pressure, are of several orders of magnitude larger than the bending energy, and usually serve 
in computational models as area and enclosed volume constraints of the bilayer membrane, 
respectively. With this simplification, often referred to as the spontaneous curvature model, the deformation energy 
of the bilayer membrane is given by the bending curvature equation with 
two constraints. The equilibrium configuration of lipid bilayer vesicles can be obtained by minimizing this 
energy. For sufficiently simplified membrane systems such as isolated vesicles with symmetric lipid composition, 
analytical analysis of these energies can give an excellent classification of the phases of the vesicle configurations, 
particularly the axisymmetric configurations, and can accurately locate the conditions under which the phase transition 
occurs \cite{Seifert1991,Seifert1997,OuyangZ1987a}. Giving the complexity the realistic interacted biological system with which 
lipid bilayer is interacted, analytical approaches may fail in quantifying these interactions and thus computational methods 
become indispensable.

The solutions of the energy minimization problem have been computed in various ways. Solving the Euler-Lagrange 
equations directly has so far been restricted by the surface parameterization to solutions which are axisymmetric \cite{Seifert1991, Seifert1997, Feng2006}. 
One alternative is to minimize the energy over a smaller subspace of membrane configurations. On triangulated
vesicle surface, the subspace can be that spanned by the basis functions at triangular vertices. In this
subspace the deformation energy can be approximated using Rayleigh-Ritz procedure \cite{Canham1970,HeinrichV1999a} 
or finite element methods \cite{Feng2006, Ma2008}. In case that the triangular mesh needs to be locally refined to resolve a 
very large local curvature, the number of basis functions becomes large, leading to a significant increase of 
computational cost of numerical minimization. The other alternative is to define the surface of a lipid bilayer
vesicle as the level set of a phase field function. Geometrical properties of the vesicle surface can be represented
using the phase field function, and the equilibrium surface configuration can be obtained by following the gradient 
flow driven by the membrane curvature. Approaches of this type have been very successful in describing the 
nonaxisymmetric equilibrium configurations \cite{Du2004,DuQ2005b,Du2006}, membranes with multiple lipid species \cite{WangX2008a}, membranes with surfactant 
sedimentation or solvation \cite{XuJ2012a,TeigenK2011a}, and membrane-flow interactions \cite{SohnJ2010a,LiS2012a}. Moreover, the phase field approach 
is the only computational method that can simulate topological changes of the surface during the merging or separation of vesicles. 
Nevertheless, there is a significant increase in the computational cost
with phase field approaches, because the tracking of the 2-D vesicle surface is replaced by the evolution of 3-D phase 
field function. Special numerical techniques, such as spectral methods \cite{Du2004,DuQ2005b,Du2006}, finite difference or finite element methods 
with local mesh refinement \cite{YangX2006a,WiseS2007a}, have been developed to accelerate the simulations.

There are many important applications where fast algorithms are necessary for simulating vesicle deformation that does not 
necessarily involve topological change. These include the vesicle deformations induced by protein clusters \cite{Farsad2003,SolonJ2005a}, and
interactions between blood flow in micro-vessels and the deformable RBCs \cite{EggletonC1998a,BagchiP2007a}. For example, about 
7778 nodes are used to describe the deformation of spherical vesicles in 3-D simulations of single vesicle-blood flow interactions using
immersed boundary method \cite{EggletonC1998a}. In a 2-D simulation of multiple vesicle-blood flow interactions, about 128 to 512 nodes
are used to track the deformation of a single vesicle at a satisfactory resolution \cite{BagchiP2007a}. A full 3-D simulation of this type
can be even more expensive. An accurate and simplified representation of vesicle configuration will greatly improve the 
efficiency of these computational investigations.

In this article, we approximate the deformation energy of lipid vesicles using surface harmonic functions. These real-valued 
functions are the linear combination of complex-valued spherical harmonics functions. They share the same orthonormal and completeness 
properties as spherical harmonics, and thus provide a natural basis for vesicle configurations of star-shape. 
The uniform convergence of spherical harmonic functions in approximating smooth functions 
enables us to choose far fewer number of basis 
functions to approximate vesicle surface when compared to finite element methods with a triangular mesh. With the surface harmonic expansion, 
one can analytically compute the variations of the deformation energy with respect to the expansion coefficients. 
A nonlinear conjugate gradient method is employed to numerically compute the minimizer of the deformation energy, without incurring the computationally cost evaluation of Hessian 
matrix for the deformation energy. Our numerical experiments below show that 49 surface harmonic functions are sufficient to represent 
axisymmetric prolate, oblate, and stomatocyte shapes as well as many nonaxisymmetric shapes. We note that complex-valued spherical harmonic 
functions were used to compute the local curvature of vesicles deformed by protein clusters, where least square solution of an over-determined
linear systems for spherical harmonic expansion coefficients is solved using singular value decomposition \cite{BahramiA2011a}. In a more relevant study, a 
complex-valued spherical harmonic expansion is sought for each individual Cartesian coordinate of the node on vesicle surface \cite{KhairyK2011a}.

Curvature energy by itself can not account for the deformation and various biological functions of the membrane, 
the potential and kinetic energies arising in realistic intermolecular interactions shall be included in simulations \cite{DasS2008a,Mikucki2014}. 
Although the vesicle deformation can be well represented using spherical harmonic functions as shown in this study, whether the 
modeling of various intermolecular interactions and the total energy can be facilitated by the spherical harmonic expansions is to be determined. 
For fluid flows, a simplified geometric representation of immersed vesicles will speed up the evaluation of flow field. 
For Stokes flow in particular, the implementation of the Stokeslet can be accelerated by using the analytical expression of the
vesicle surfaces. In a companion paper we will show that our fast algorithm can be applied to minimize a total energy
that ensembles both the deformation energy and the electrostatic potential energy arising from protein-membrane
interactions, where a transfer matrix assigns the electrostatic force computed at an arbitrary point 
on vesicle using surface harmonic functions.

The rest of the article is organized as follows.  In Section 2, we present the bending energy of a lipid bilayer membrane and 
its variation along with the area and volume constraints. Section 3 introduces surface harmonics and establishes the 
variational problem in terms of the surface harmonics parameterization. Section 4 presents examples of equilibrium 
configurations and confirms these results from previous work in the literature. Finally, Section 5 discusses future developments 
of the model, including electrostatics from protein-membrane interactions.  


\section{Lipid bilayer mechanics}

\subsection{Total energy and constraints}

According to the classic spontaneous curvature model as developed by Canham \cite{Canham1970}, Helfrich \cite{Helfrich1973}, and Evans \cite{Evans1974}, the bending energy of the membrane is given by 
\begin{equation}\label{eq:bending_energy}
E[\Gamma] = \int_{\Gamma} \left(\frac{1}{2} \mathscr{K}_C (2H-C_0)^2 +\mathscr{K}_G K \right)  dS, 
\end{equation}
where $\mathscr{K}_C$ and $\mathscr{K}_G$ are the bending modulus and Gaussian saddle-splay modulus, respectively, $H$ is the mean curvature, $C_0$ is the spontaneous curvature, and $K$ is the Gaussian curvature \cite{Feng2006}.  The position of the membrane is given by $\Gamma$, and so the total bending energy $E[\Gamma]$ depends only upon the membrane position.  An explicit formula for the differential surface area element $dS$ will be given in terms of the surface harmonic parameterization in Section 3.  We assume that the spontaneous curvature of the membrane $C_0$ is a constant.  

In the case of protein-membrane interactions, the membrane will also be under an external force. The potential energy from this interaction may be added to the total potential energy of the system, denoted by $\Pi[\Gamma]$:
\begin{equation}\label{eq:II}
\Pi[\Gamma] = E[\Gamma] + G[\Gamma],
\end{equation}
where $E[\Gamma]$ is the bending energy \eqref{eq:bending_energy} and $G[\Gamma]$ is the electrostatic potential energy from the protein-membrane interaction.  The introduction of the electrostatic potential energy will be discussed in the companion paper.  For this paper, we only consider the bending energy $E[\Gamma]$.  

To complete the total energy and account for the area expansion/contraction and osmotic pressure, the constraints for the conservation of the surface area of the membrane and the total volume enclosed are added to \eqref{eq:II}.  The equilibrium position of the bilayer membrane is determined by the surfaces which minimize the bending energy with the surface area and volume constraints.  The total potential energy (neglecting any electrostatic potential energy for now) with the penalty functions is given by 
\begin{equation}\label{eq:I}
I[\Gamma] = E[\Gamma] + \frac{k_S}{2}(S_A-\bar{S})^2 + \frac{k_V}{2}(V-\bar{V})^2
\end{equation}
where $S_A$ and $V$ are the total surface area and volume, $\bar{S}$ and $\bar{V}$ are the initial surface area and volume, and $k_S$ and $k_V$ are large constants to enforce the constraints to a chosen degree of precision.

\subsection{Variational formulation}

The minimization of \eqref{eq:I} gives rise to the bending curvature equation of $\Gamma$, 
\begin{equation}\label{eq:dI}
\delta_\Gamma I[\Gamma] = \delta_\Gamma E[\Gamma] + k_S(S_A-\bar{S})\delta S_A + k_V(V-\bar{V})\delta V.
\end{equation}
The computation of the terms in \eqref{eq:bending_energy}, \eqref{eq:I}, and \eqref{eq:dI} depend on the choice of parameterization of the membrane surface.  A convenient consequence of the bending energy of lipid bilayers is that the energy functionals are independent of the surface parameterization \cite{Capovilla2003}.  Therefore, we choose to parameterize the surface not by brute force, intrinsic cartesian coordinates, but by surface harmonics.  Since lipid vesicles are sphere-like structures, the choice of surface harmonics to represent the surface is natural, and this choice reduces the number of terms necessary for the computation.  


\section{Surface harmonics approximation}

In this section, we introduce the surface harmonic parameterization, and provide formulas for the terms in \eqref{eq:bending_energy}, \eqref{eq:I}, and \eqref{eq:dI} according to this parameterization.  Surface harmonic functions are a simplification of spherical harmonic functions.  We choose to implement surface harmonics in our work to allow for an easier minimization algorithm.  Spherical harmonic functions are linear combinations of real valued spherical harmonics. 

\subsection{Spherical harmonic functions}

Spherical harmonics are solutions to Laplace's Equation in spherical coordinates.  The solution can be obtained through separation of the variables $\theta$ and $\phi$; however, more convenient way to construct spherical harmonics is to use a generalization of Legendre polynomials.  Legendre polynomials, also called Legendre functions of the first kind, Legendre coefficients, or zonal harmonics, are solutions to the Legendre differential equation.  The Legendre polynomial can be defined by the contour integral
\begin{equation}\label{eq:Lpoly}
P_n(z) = \frac{1}{2\pi i}\oint (1-2tz+t^2)^{-1/2}t^{-n-1} \; dt.
\end{equation}
Another useful representation utilizes the Rodrigues representation,
\begin{equation}\label{eq:rodrigues}
P_n(x) = \frac{1}{2^n n!} \frac{d^n}{dx^n}(x^2 - 1)^n.
\end{equation}

The associated Legendre polynomials generalize Legendre polynomials, provided $m\neq0$, and are defined by
\begin{equation}\label{eq:ALpoly}
P^m_n(x) = (-1)^m(1-x^2)^{m/2} \frac{d^m}{dx^m}P_n(x),  \qquad m\neq 0.
\end{equation}
If $m=0$, the associated Legendre polynomial is just the Legendre polynomial.  By Rodrigues' formula, 
\begin{equation}\label{eq:ALpolyR}
P^m_n(x) = (-1)^m(1-x^2)^{m/2} \frac{d^m}{dx^m}\left( \frac{1}{2^n n!} \frac{d^n}{dx^n}(x^2 - 1)^n\right), \qquad m\neq 0
\end{equation}
It is convenient to introduce the change of variables $\mu=\cos(\theta)$.  In this way, the partial derivatives with respect to the polar angle $\theta \in [0,\pi]$ may be computed.  Using this notation, normalized spherical harmonic functions are defined by
\begin{equation}\label{eq:SHF}
Y_n^m(\theta,\phi) = \sqrt{\left( \frac{(2n+1)(n-m)!}{4\pi(n+m)!}\right)} P^m_n(\mu) e^{im\phi},
\end{equation}
where $P^m_n(\mu)$ is the associated Legendre polynomial evaluated at $\mu = \cos(\theta)$.  
%

Since spherical harmonics form an orthonormal basis for $L^2(\reals^2)$, linear combinations of them can represent smooth surfaces.  The surface is parameterized in spherical coordinates $(\theta, \phi, r(\theta,\phi))$, where the radius $r$ is expressed in terms of spherical harmonics,
\begin{equation}\label{eq:r_SHF}
r(\theta,\phi) = \sum_{n=0}^\infty \sum_{m=-n}^n a_n^m Y^m_n(\theta,\phi).
\end{equation}
The $a_n^m$ are the coefficients of the linear representation.  These coefficients can be determined by the following formula:
\begin{equation}\label{eq:anm_SHF}
a_n^m = \int_0^{2\pi} \int_0^\pi r(\theta,\phi)\overline{Y^m_n(\theta,\phi)}\sin(\theta) \; d\theta\, d\phi
\end{equation}
where $\overline{Y^m_n(\theta,\phi)}$ is the complex conjugate of $Y^m_n(\theta,\phi)$.  

\subsection{Surface harmonic functions}

If the coefficients $a_n^m$ are poorly chosen so that the complex parts of $a_n^m$ and $Y_n^m$ do not cancel the radius parameterizing the object \eqref{eq:r_SHF} will be complex.  In an optimization routine, the coefficients are perturbed arbitrarily, so any nonzero perturbation in the complex part will result in a complex surface.  Since we seek a real-valued surface that minimizes the potential energy \eqref{eq:I}, we use only the real parts of the spherical harmonics to ensure that the surface under the energy minimization is real.  These are surface harmonics. 


%

By Euler's formula, each spherical harmonic function can be rewritten as 
\begin{equation}\label{eq:Y_Eulers}
Y^m_n(\theta,\phi) =  f^m_n P^m_n(\mu)(\cos(m\phi) + i\sin(m\phi)).
\end{equation}
where $f^m_n$ is the normalization factor
\begin{equation}\label{eq:fnm}
f^m_n = \sqrt{\left( \frac{(2n+1)(n-m)!}{4\pi(n+m)!}\right)}.
\end{equation}
Since the real and complex parts of \eqref{eq:Y_Eulers} are solutions to Laplace's equation, we 
define the surface harmonics to be
\begin{equation}\label{eq:Snm}
S_n^m(\theta, \phi) = \left\{\begin{split} 
f_n^m P^{m}_n(\mu)\cos(m\phi) & \qquad \textrm{if } m\geq0 \\
f_n^{|m|} P^{|m|}_n(\mu)\sin(|m|\phi) & \qquad \textrm{if } m <  0 
\end{split}\right..
\end{equation}
%
We now parameterize the radius of a smooth surface by a linear combination of the surface harmonics $S_n^m(\theta,\phi)$, 
\begin{equation}\label{eq:r}
r(\theta,\phi) = \sum_{n=0}^\infty \sum_{m=-n}^n a_n^m S_n^m(\theta, \phi).
\end{equation}
%
A more convenient representation of the radius of a surface that avoids the sign changes in $m$ is
\begin{equation}\label{eq:r1}
r(A_n^m, B_n^m; \theta,\phi) = \sum_{n=0}^\infty \sum_{m=0}^n \Big(A_n^m \cos(m\phi) + B_n^m \sin(m\phi)\Big)f_n^mP_n^m(\mu),
\end{equation}%
where
\begin{equation}\label{eq:anm}
a_n^m = \left\{\begin{split}
A_n^m & \qquad \textrm{if } m\geq0\\
B_n^m & \qquad \textrm{if } m<0
\end{split}\right..
\end{equation}
%
%

%


\subsection{Discretizing the surface}


Discretize the surface by $n_t$ values of $\theta$ and $n_p$ values of $\phi$ for a total of $\mathcal{N} =n_t\cdot n_p$ points.  Since $\mathcal{N}$ can be a very large number, performing pointwise calculations on the mesh can be computationally costly.  Under the surface harmonics parameterization, the surface is approximated by truncating the infinite sum in \eqref{eq:r} at some number $N$.  This reduces the computation cost since there are far fewer surface harmonic functions required to approximate a surface than a curvilinear cartesian grid of mesh points.  For fixed values of $\theta_k$, $k=1,\cdots,n_t$, and $\phi_l$, $l=1,\cdots,n_p$, the radius $r$ is determined by the coefficients $A_n^m$ and $B_n^m$, 
\begin{equation}\label{eq:r_trunc}\begin{split}
r_{kl}(A_n^m, B_n^m) &= r_{kl}(A_n^m, B_n^m; \theta_k, \phi_l) \\
                     &= \sum_{n=0}^N \sum_{m=0}^n \Big(A_n^m \cos(m\phi_l) + B_n^m \sin(m\phi_l)\Big)f_n^mP_n^m(\mu).
\end{split}\end{equation}
Using this parameterization, there are 
a total of $(N+1)^2$ coefficients to determine the surface parameterization.  Our numerical results confirm that $(N+1)^2 << \mathcal{N}$.
For notational convenience, let $\vec{a}$ be a vector of all of the surface harmonic coefficients given by \eqref{eq:anm}.  
\begin{equation}\label{eq:a}
\vec{a} = [A_0^0, \; A_1^0, \;A_1^1, \;B_1^1, \;A_2^0, \;A_2^1, \; A_2^2, \; B_2^1, \; B_2^2, \; \cdots, \; A_N^0, \; \cdots, \; A_N^N, \; B_N^1, \; \cdots, \; B_N^N]^T
\end{equation}
We use the subscript $i = 0, \cdots, (N+1)^2-1$ to denote the surface harmonic mode $a_i$.  
Next, we use the surface harmonics parameterization \eqref{eq:r_trunc} to finish the formulas for \eqref{eq:bending_energy}, \eqref{eq:I}, and \eqref{eq:dI}.  

\subsection{Energy formulation in terms of surface harmonic parameterization}

The differential surface element 
is 
\begin{equation}\label{eq:SHF_S}
dS = r [r_{\phi}^2+r_{\theta}^2 \sin^2(\theta) +r^2 \sin^2(\theta)]^{1/2} \;d\theta \,d\phi ,
\end{equation}
and so the surface area of $\Gamma$ is 
\begin{equation}\label{eq:SA1}
S_A = \int_0^{2\pi} \int_0^{\pi} r [r_{\phi}^2+r_{\theta}^2 \sin^2(\theta) +r^2 \sin^2(\theta)]^{1/2} d\theta d\phi .
\end{equation}
For simplicity in later calculations, we define the determinant of the covariant metric tensor to be 
\begin{equation}\label{eq:omega}
\omega = r [r_{\phi}^2+r_{\theta}^2 \sin^2(\theta) +r^2 \sin^2(\theta)]^{1/2},
\end{equation}
so that $dS = \omega \, d\theta d\phi$, and  
\begin{equation}\label{eq:SA}
S_A = \int_0^{2\pi} \int_0^{\pi} \omega \;d\theta d\phi .
\end{equation}
The volume enclosed by the membrane is 
%
%
\begin{equation}\label{eq:Vol}
V = \frac{1}{3}\int_0^{2\pi} \int_0^\pi r^3\sin(\theta)\;d\theta\;d\phi .
\end{equation}
These two formulas for surface area and volume parameterize the constraints in \eqref{eq:I}.  The variation of \eqref{eq:SA} with respect to a surface harmonic mode $a_i$ can be computed directly as
\begin{equation}
\delta S_A = \int_0^{2\pi} \int_0^\pi \frac{\delta \omega}{\delta a_i} \; d\theta \; d\phi.
\end{equation}
Similarly, the variation of \eqref{eq:Vol} can be computed directly,  
\begin{equation}
\delta V = \int_0^{2\pi} \int_0^\pi r^2 \frac{\delta r}{\delta a_i} \sin(\theta) \; d\theta \; d\phi.
\end{equation}
%

The derivatives of $r$ given by \eqref{eq:r_trunc} with respect to $\theta$ and $\phi$ are computed next.  The derivatives with respect to $\phi$ are straightforward.
The subscripts in the following formulas represent partial derivatives and should not be confused with mesh positions $k$ and $l$.  
\begin{align}
r_{\phi} &= \sum_{n=0}^N \sum_{m=0}^n \Big( -m\, A_n^m\sin(m\phi) + m\, B_n^m\cos(m\phi) \Big) f_n^m P^m_n(\mu) \notag\\
         &= \sum_{n=0}^N \sum_{m=-n}^{n} -m\, a_n^{m} S_n^{-m}\label{eq:r_phi}\\
r_{\phi\phi} &= \sum_{n=0}^N \sum_{m=0}^n \Big(-m^2 A_n^m\cos(m\phi) + -m^2 B_n^m\sin(m\phi) \Big) f_n^m P^m_n(\mu) \notag\\
             &=-m^2r\label{eq:r_phiphi}
\end{align}
The derivatives with respect to $\theta$ are
\begin{equation}\label{eq:r_theta}
r_{\theta} = \sum_{n=0}^N \sum_{m=0}^n \Big( A_n^m\cos(m\phi) + B_n^m\sin(m\phi) \Big) f_n^m \partial_\theta P_n^m(\mu)
\end{equation}
\begin{equation}\label{eq:r_thetatheta}
r_{\theta \theta} = \sum_{n=0}^N \sum_{m=0}^n \Big( A_n^m\cos(m\phi) + B_n^m\sin(m\phi) \Big) f_n^m \partial_\theta^2 P_n^m(\mu) 
\end{equation}
\begin{align}
r_{\phi\theta} = r_{\theta\phi} &= \sum_{n=0}^N \sum_{m=0}^n \Big( -m \,A_n^m\sin(m\phi) + m\, B_n^m\cos(m\phi) \Big) f_n^m \partial_\theta P_n^m(\mu) \notag\\
                                &= \sum_{n=0}^N \sum_{m=-n}^{n} -m\, a_n^{m} \frac{\partial S_n^{-m}}{\partial \theta}, \label{eq:r_phitheta}
\end{align}
where $\partial_\theta P_n^m(\mu)$ is given by the recurrence relation for the derivative of the associated Legendre polynomial $P_n^m(\mu)$,
\begin{equation}\label{eq:dAL}
\partial_{\theta} P^m_n(\mu) = \frac{-1}{sin(\theta)}\left( (n+1) \cos(\theta) P^m_n(\mu) - (n-m+1) P^m_{n+1}(\mu) \right),
\end{equation}
and $\partial_\theta^2 P_n^m(\mu)$ is computed directly from \eqref{eq:dAL} as 
\begin{equation}\label{eq:d2AL}
\begin{split}
\partial_\theta^2 P_n^m(\mu) &= \Big((n+1+(n+1)^2\cos^2{\theta}) P_n^m(\mu) - 2\cos{\theta}(n-m+1)(n+2)P_{n+1}^m(\mu) \\
& \qquad   + (n-m+1)(n-m+2)P^m_{n+2}(\mu)\Big)\frac{1}{\sin^2{\theta}}.
\end{split}
\end{equation}
%

The variation of $r$ is given by 
\begin{equation}\label{eq:drda1}
\frac{\delta r}{\delta a_i} = \left\{\begin{split}
                              \frac{\delta r}{\delta A_n^m} &= f_n^m P_n^m(\mu)\cos(m\phi)&\qquad m \geq 0\\
                              \frac{\delta r}{\delta B_n^m} &= f_n^{|m|} P_n^{|m|}(\mu)\sin(|m|\phi)&\qquad m < 0
\end{split}\right.\end{equation}
We notice immediately that these variations match \eqref{eq:Snm} exactly, and so the variations are just surface harmonic functions.  For simplicity in the formulas, for the surface harmonic mode $a_i \in \vec{a}$, $i = 0, 1, 2, \cdots (N+1)^2-1$, we define $n$ and $m$ from $i$ to match \eqref{eq:Snm} and \eqref{eq:drda1} as follows:
\begin{equation}\label{eq:compute_n}
n(i) = \lfloor\sqrt{i}\rfloor 
\end{equation}
\begin{equation}\label{eq:compute_m}
m(i) = \left\{\begin{split} &i-n^2   & \qquad &\textrm{if }(i-n^2) \leq n \\
                         &n^2+n-i & \qquad &\textrm{otherwise}
     \end{split}\right.
\end{equation}
Then, the variation in $r$ with respect to the mode $a_i$ is given by 
\begin{equation}\label{eq:drda}
\frac{\delta r}{\delta a_i} = S_n^m,
\end{equation}
which through \eqref{eq:compute_n} and \eqref{eq:compute_m}, only depends on the mode $i$.  
Continuing with this simplification of subscripts, 
\begin{equation}\label{eq:dr_phida}
\frac{\delta r_{\phi}}{\delta a_i} = -m S_n^{-m} 
\end{equation}
\begin{equation}\label{eq:dr_phiphida}
\frac{\delta r_{\phi\phi}}{\delta a_i} = -m^2 \frac{\delta S}{\delta a_i}  
\end{equation}
The variations of $r_\theta$ and $r_{\theta\theta}$ are 
\begin{equation}\label{eq:dr_thetada}
\frac{\delta r_\theta}{\delta a_i} = 
\left\{\begin{split}
\frac{\delta r_\theta}{\delta A_n^m} &= \cos(m\phi)f_n^m \partial_\theta P_n^m(\mu) &\qquad m \geq 0\\
\frac{\delta r_\theta}{\delta B_n^m} &= \sin(|m|\phi)f_n^{|m|}  \partial_\theta P_n^{|m|}(\mu)&\qquad m < 0
\end{split}\right.\end{equation}
\begin{equation}\label{eq:dr_thetathetada} 
\frac{\delta r_{\theta\theta}}{\delta a_i} = 
\left\{\begin{split} 
\frac{\delta r_{\theta\theta}}{\delta A_n^m} &= \cos(m\phi)f_n^m \partial_\theta^2 P_n^m(\mu) &\qquad m\geq 0\\
\frac{\delta r_{\theta\theta}}{\delta B_n^m} &= \sin(|m|\phi)f_n^{|m|}  \partial_\theta^2 P_n^{|m|}(\mu) & \qquad m < 0
\end{split}\right.\end{equation}
\begin{equation}\label{eq:dr_phithetada}
\frac{\delta r_{\theta\phi}}{\delta a_i} = 
\left\{\begin{split}
\frac{\delta r_{\theta\phi}}{\delta A_n^m} &= \frac{\delta r_{\phi\theta}}{\delta A_n^m} = -m \sin(m\phi)f_n^m \partial_\theta P_n^m(\mu) & \qquad m \geq 0\\
\frac{\delta r_{\theta\phi}}{\delta B_n^m} &= \frac{\delta r_{\phi\theta}}{\delta B_n^m} = -m \cos(|m|\phi)f_n^{|m|}  \partial_\theta P_n^{|m|}(\mu)& \qquad m < 0
\end{split}\right.\end{equation}
Finally, 
the variation of $\omega$ is 
%
\begin{equation}\label{eq:dw}
\frac{\delta \omega}{\delta a_i} = 
\left\{\begin{split}
\frac{\delta \omega}{\delta A_n^m} &= \frac{\delta r}{\delta A_n^m}[r_{\phi}^2+r_{\theta}^2 \sin^2(\theta) +r^2 \sin^2(\theta)]^{1/2} \\
	& \qquad + \frac{r}{2} \left(\frac{ 2r_\phi \DD\frac{\delta r_\phi}{\delta A_n^m} +2r_{\theta}\frac{\delta r_\theta}{\delta A_n^m} \sin^2(\theta) +2r \frac{\delta r}{\delta A_n^m}\sin^2(\theta) }{\sqrt{r_{\phi}^2+r_{\theta}^2 \sin^2(\theta) +r^2 \sin^2(\theta)}}\right) \qquad m \geq 0\\
\frac{\delta \omega}{\delta B_n^m} &= \frac{\delta r}{\delta B_n^m}[r_{\phi}^2+r_{\theta}^2 \sin^2(\theta) +r^2 \sin^2(\theta)]^{1/2} \\
	& \qquad + \frac{r}{2} \left(\frac{ 2r_\phi \DD\frac{\delta r_\phi}{\delta B_n^m} +2r_{\theta}\frac{\delta r_\theta}{\delta B_n^m} \sin^2(\theta) +2r \frac{\delta r}{\delta B_n^m}\sin^2(\theta) }{\sqrt{r_{\phi}^2+r_{\theta}^2 \sin^2(\theta) +r^2 \sin^2(\theta)}}\right)\qquad m < 0
\end{split}\right.\end{equation}


To finish the formulas, we need expressions for the mean and Gaussian curvatures $H$ and $K$, respectively, and their variations, in terms of $r$.  It is convenient to first define the so-called ``warping functions" $E, F,$ and $G$, and shape operator functions $L, M$, and $N$.  The warping functions are the coefficients of the first fundamental form and are given by
\begin{align}
E&=r^2_{\theta}+r^2 \label{eq:E}\\
F&=r_{\theta} r_{\phi} \label{eq:F}\\
G&=r^2_{\phi}+r^2 \sin^2(\theta) \label{eq:G}
\end{align}
The shape operator functions are the coefficients of the second fundamental form and are given by
\begin{align}
L&=-\vec{x}_{\theta} \cdot \hat{n}_{\theta} \label{eq:Ldef} \\
M&=\frac{1}{2} \left (\vec{x}_{\theta} \cdot \hat{n}_{\phi} + \vec{x}_{\phi} \cdot \hat{n}_{\theta} \right)  \label{eq:Mdef}\\
N&=-\vec{x}_{\phi} \cdot \hat{n}_{\phi}\label{eq:Ndef}
\end{align}
%
where $\hat{n}$ is the unit normal to the surface at $\vec{x}$,
\begin{equation}\label{eq:SHF_n}
\hat{n} = \frac{\vec{x}_\theta \times \vec{x}_\phi}{\left|\vec{x}_\theta \times \vec{x}_\phi\right|}.
\end{equation}
Define 
\begin{equation}\label{eq:R}
R= \left|\vec{x}_\theta \times \vec{x}_\phi\right|
\end{equation}
for notational convenience.  
The shape operator functions can be expressed in terms of $r$ and its derivatives.  
\begin{align}
L&= R^{-1} (-2r r_\theta^2\sin(\theta)+r^2 r_{\theta\theta}\sin(\theta) - r^3\sin(\theta)) \label{eq:L} \\
M&= R^{-1} (2r r_\phi r_\theta \sin(\theta) - r^2r_{\theta\phi} \sin(\theta) + r^2 r_\phi \cos(\theta)) \label{eq:M}\\
N&= R^{-1} (-r^3\sin^3(\theta) + r^2r_{\phi\phi} \sin(\theta) + r^2r_\theta \cos(\theta)\sin^2(\theta) - 2r r_\phi^2\sin(\theta)) \label{eq:N}
\end{align}
Finally, we obtain the local mean curvature and the Gaussian curvature,
\begin{align}
H(\theta,\phi) &= \frac{EN+GL-2FM}{2(EG-F^2)} \label{eq:H}\\
K(\theta,\phi) &= \frac{LN-M^2}{EG-F^2} \label{eq:K}
\end{align}
In terms of $r$, 
\begin{equation}\label{eq:H_expr}\begin{split}
H(\theta,\phi) &= \frac{1}{2}R^{-1} \left[ 3r_\theta^2 r^2 \sin^3(\theta) - r_\theta^2 r r_{\phi\phi} \sin(\theta) -r_\theta^3 r \cos(\theta)\sin^2(\theta) + 8 r_\theta^2 r_\phi^2 \sin(\theta) \right.\\
	& \phantom{\frac{1}{2}}\qquad + 2r^4\sin^3(\theta) - r^3 r_{\phi\phi} \sin(\theta) -r^3 r_\theta \cos(\theta)\sin^2(\theta) + 3r^2 r_\phi^2\sin(\theta) \\
	&\phantom{\frac{1}{2}}\qquad \left.- r_\phi^2 r r_{\theta\theta} \sin(\theta) - r^3 r_{\theta\theta}\sin^3(\theta) - 2 r_\phi r_\theta r r_{\theta\phi} \sin(\theta) + 2r_\phi^2 r_\theta r \cos(\theta) \right] \\
	& \phantom{\frac{1}{2}}\qquad / \left[ -r^3\sin^2(\theta) - r r_\theta^2\sin^2(\theta) - r r_\phi^2 \right].
\end{split}\end{equation}
The variation of the mean curvature with resepect to a suface harmonic mode $a_i$ is striaghtforward.

Finally, we compute the variation of $E[\Gamma]$.  To do this, we assume that the Gaussian modulus $\mathcal{K}_G$ is uniform over the membrane surface, and so the Gaussian curvature integrates to a constant $\int_\Gamma K dS = 4\pi(1-g)$ where $g$ is the genus of the membrane topology, according to the Gauss-Bonet Theorem \cite{Sokolnikoff1964}.  Thus, the variation of the bending energy with respect to a surface harmonic mode $a_i$ is
\begin{equation}\label{eq:d_bending_energy}\begin{split}
\delta_\Gamma E[\Gamma; \theta,\phi] = \int_\Gamma \left[ \mathscr{K}_C 2(2H-C_0) \delta H \, \omega + \mathscr{K}_C\frac{1}{2}(2H-C_0)^2\delta\omega \right] \,d\theta d\phi.
\end{split}\end{equation}

The variational form of the total energy with respect to spherical harmonic coefficients, given by \eqref{eq:dI}, is now complete.  With the surface $\Gamma$ expressed in terms of the surface harmonic coefficients $\vec{a}$, the bending curvature equation is 
\begin{equation}\label{eq:d_pi}
\delta_\Gamma I[\Gamma(\vec{a})] = \delta_\Gamma E[\Gamma(\vec{a})] + k_S(S_A-\bar{S})\delta S_A + k_V(V-\bar{V})\delta V.
\end{equation}

\subsection{Numerical methods}
We employ a Fletcher-Reeves type nonlinear conjugate gradient (NCG) method to minimize the total 
]energy functional \eqref{eq:I}.  For the parameter $\beta$, we chose the Hestenes-Stiefel formula.  For a 
description of the method, please refer to \cite{NumericalOptimization}. The psuedocode is provided in 
the Appendix \ref{append}.
 
\subsection{Expansion modes and quadrature}

In practice, the surface harmonic expansion \eqref{eq:r_trunc}  is truncated at some value $N$, giving a total of $(N+1)^2$ surface harmonic modes used.  To determine an appropriate $N$, we reconstructed three surfaces and examined the root mean square error in the surface reconstruction pointwise, and the relative error in the volume, surface area, and energy.  For many vesicle structures, increasing $N$ ahieves higher accuracy; however, it also increases computation time.  For all structures, the accuracy is perfect as $N\to\infty$, however; the convergence is not monotone.  For some vesicle structures, increasing $N$ transiently actually gave a worse approximation.  For these reasons, we chose to perform the energy minimization procedure with the smallest possible $N$ giving a desired accuracy.  

\begin{figure}[h!tbs]
    \centering
   \includegraphics[height=3cm]{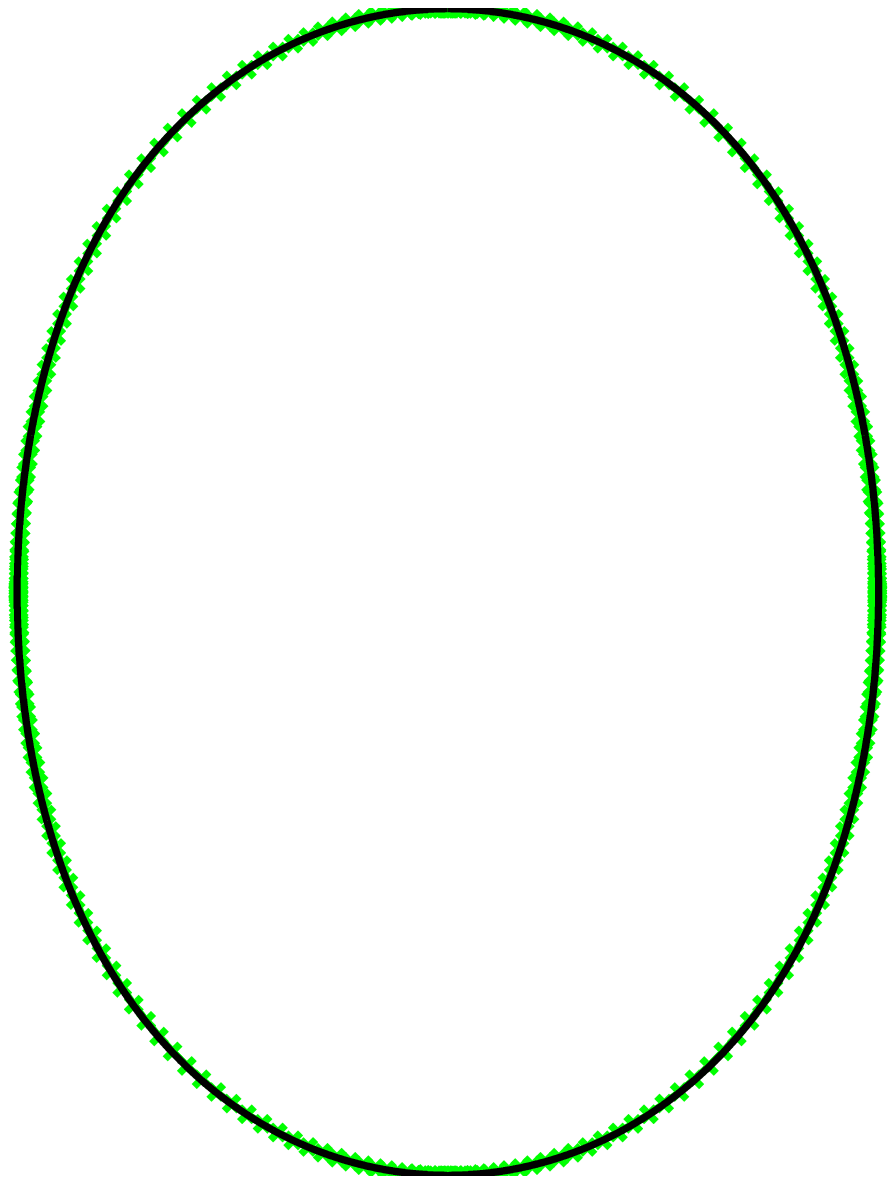} 
   \includegraphics[height=3cm]{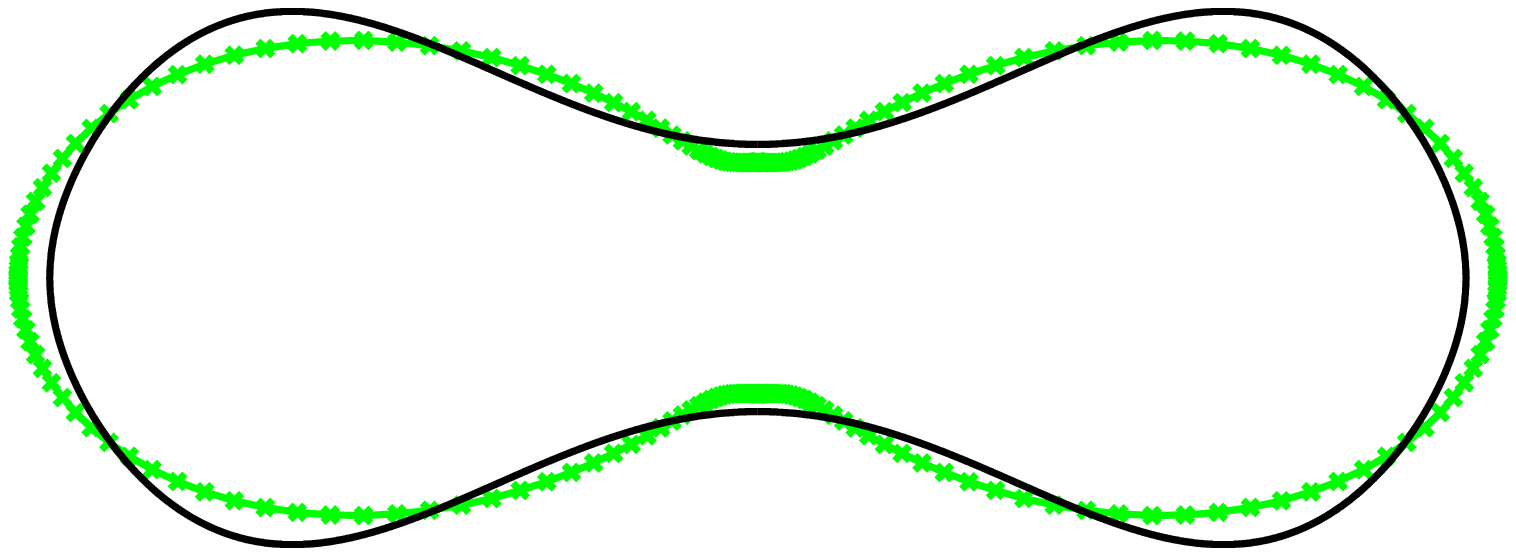} \qquad
   \includegraphics[height=3cm]{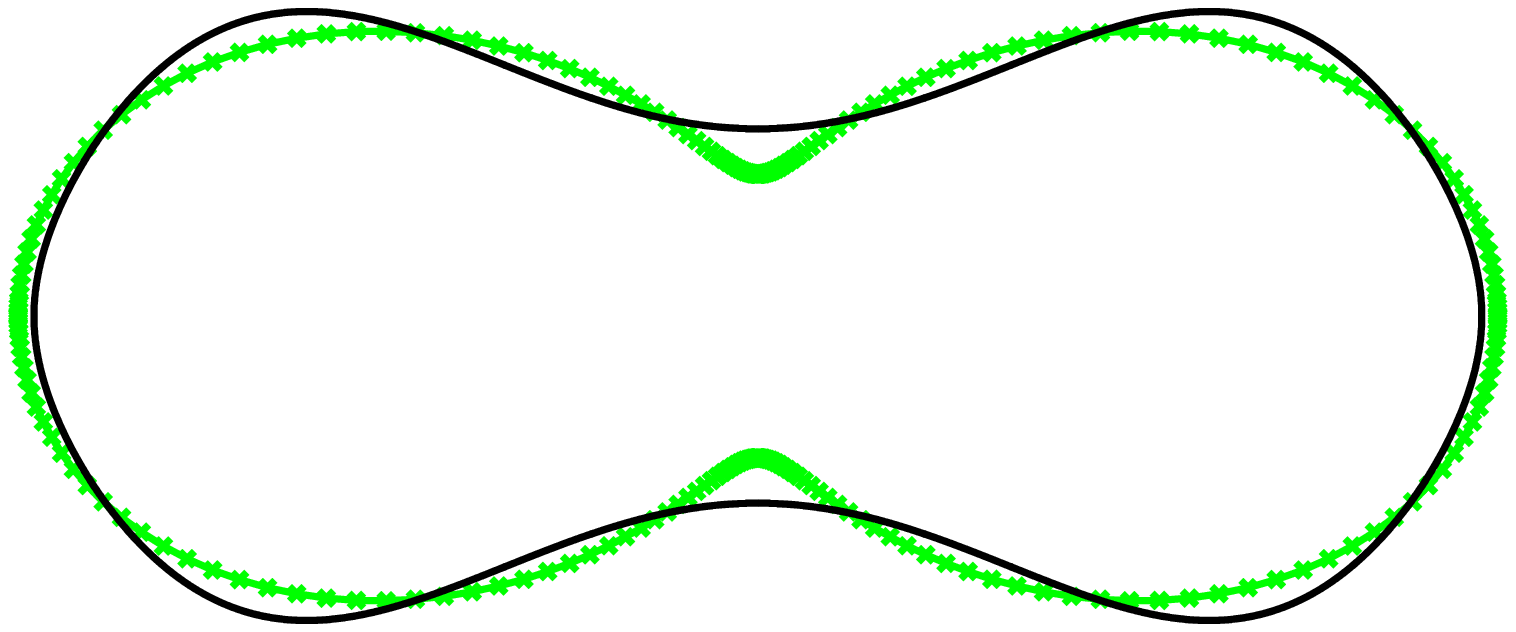} 
\caption{Profile of various test surfaces (black) and their reconstructions using surface harmonics with $N=4$ (green).  The full surfaces are generated by rotation about thet $y$-axis.  Left: prolate from \cite{Seifert1991}, middle and right: red blood cell from \cite{Evans1972} with 50\% and 90\% weight on tonicity 217 mO coefficients, respectively.}  \label{fig:samples} 
\end{figure}

The first surface we reconstructed was an energy minimizing axisymmetric prolate from Seifert et. al., \cite{Seifert1991}.  Instructions for reconstructing this surface can be found in Appendix B of \cite{Seifert1991}, with choice of parameters $\bar{P} = 0.1$, $\bar{\Sigma} = -1.1\bar{P}^{2/3}$, $C_0 = 0$, and $U(0) = 0.56$.  
Next, we reconstructed statistically fitted parameterizations of a red blood cell from \cite{Evans1972}.  The height of the profile of the surface is given by
\begin{equation}\label{eq:Fung}
h(x) = \frac{\pm0.5}{R_0}(1-x^2)(C_0 + C_2x^2 + C_4x^4) \qquad x \in [-1,1].
\end{equation}
Table 4 in \cite{Evans1972} includes values for $R_0, C_0, C_2,$ and $C_4$ for producing red blood cell shapes with tonicities 300 and 217 mO.  The values are reproduced in Table \ref{tab:Evans}.  
\begin{table}[h,t,b]
\begin{center}
\begin{tabular}{lllll} \hline
Tonicity (mO) & $R_0$ ($\mu$m) & $C_0$ ($\mu$m) & $C_2$ ($\mu$m) & $C_4$ ($\mu$m)   \\
\hline
300 & 3.91 & 0.81 & 7.83 & -4.39 \\
217 & 3.80 & 2.10 & 7.58 & -5.59 \\
\hline
\end{tabular}
\end{center}
\caption{Shape coefficients for average RBC}
\label{tab:Evans}
\end{table}
We chose two linear combinations of the parameters given for averaged shapes from the ones in \cite{Evans1972}.  The profiles for the three sample surfaces and their reconstructed surfaces with $N=5$ are included in Figure \ref{fig:samples}.

The surface harmonic parameterizations of these three surfaces are determined by 
\begin{equation}\label{eq:anm_S}
a_n^m =\DD \int_0^{2\pi} \int_0^\pi r(\theta,\phi) S_n^m(\theta,\phi) \; d\theta \,d\phi.
\end{equation}
For the reconstruction, the integration was computed numerically over 230 cubature points.  Using the coefficients from \eqref{eq:anm_S}, the reconstructed radius $\tilde{r}$ was determined by \eqref{eq:r_trunc}.  The root mean square distance error in the reconstruction is defined over the cubature points by 
\begin{equation*}
E_{rmsd} = \DD\sum_{i=1}^{230} (r(\theta_i,\phi_i)-\tilde{r}(\theta_i,\phi_i))^2 / 230
\end{equation*}
The pointwise error and the relative error in the volume, surface area, and energy are provided in Tables \ref{tab:reconstruct1}-\ref{tab:reconstruct3} for various $N$. 
\begin{table}[h,t,b]
\begin{center}
\begin{tabular}{lllll} \hline
$N$ & $E_{rms}$ & $E_{Vol}$ & $E_{SA}$ & $E_{Eng}$   \\
\hline
1 & $2.06\times10^{-1}$ & $2.57\times10^{-2}$ & $2.22\times10^{-2}$ & $3.73\times10^{-2}$ \\
2 & $8.59\times10^{-4}$ & $2.22\times10^{-4}$ & $1.08\times10^{-2}$ & $8.37\times10^{-3}$ \\
3 & $8.59\times10^{-4}$ & $2.22\times10^{-4}$ & $1.08\times10^{-2}$ & $8.37\times10^{-3}$ \\
4 & $5.02\times10^{-7}$ & $2.56\times10^{-6}$ & $1.11\times10^{-2}$ & $7.19\times10^{-3}$ \\
5 & $5.06\times10^{-7}$ & $2.56\times10^{-6}$ & $1.11\times10^{-2}$ & $7.19\times10^{-3}$ \\
6 & $1.46\times10^{-7}$ & $2.76\times10^{-6}$ & $1.11\times10^{-2}$ & $7.18\times10^{-3}$ \\
7 & $1.53\times10^{-7}$ & $2.76\times10^{-6}$ & $1.11\times10^{-2}$ & $7.18\times10^{-3}$ \\
8 & $3.35\times10^{-8}$ & $2.76\times10^{-6}$ & $1.11\times10^{-2}$ & $7.18\times10^{-3}$ \\
\hline
\end{tabular}
\end{center}
\caption{Error from truncation in surface harmonic expansion for prolate sample surface.}
\label{tab:reconstruct1}
\end{table}

For the prolate surface, the reconstruction accuracy increases in all categories as $N$ increases.  
The most relevant observation to this work is that the error in the energy is less than 1\% using $N=2$ and greater.  
For a simple prolate structure, only 9 modes are required.  

\begin{table}[h,t,b]
\begin{center}
\begin{tabular}{lllll} \hline
$N$ & $E_{rms}$ & $E_{Vol}$ & $E_{SA}$ & $E_{Eng}$   \\
\hline
1 & $9.80\times10^{-2}$ & $4.00\times10^{-1}$ & $4.20\times10^{-1}$ & $4.21\times10^{-1}$ \\
2 & $6.65\times10^{-3}$ & $4.80\times10^{-2}$ & $4.99\times10^{-2}$ & $8.36\times10^{-1}$ \\
3 & $6.65\times10^{-3}$ & $4.80\times10^{-2}$ & $4.99\times10^{-2}$ & $8.36\times10^{-1}$ \\
4 & $2.29\times10^{-3}$ & $1.93\times10^{-3}$ & $1.72\times10^{-2}$ & $5.34\times10^{-2}$ \\
5 & $2.29\times10^{-3}$ & $1.93\times10^{-3}$ & $1.72\times10^{-2}$ & $5.34\times10^{-2}$ \\
6 & $1.39\times10^{-3}$ & $4.94\times10^{-3}$ & $8.35\times10^{-3}$ & $1.62$ \\
7 & $1.39\times10^{-3}$ & $4.94\times10^{-3}$ & $8.35\times10^{-3}$ & $1.62$ \\
8 & $4.13\times10^{-4}$ & $3.29\times10^{-5}$ & $1.44\times10^{-2}$ & $7.97\times10^{-1}$ \\
\hline
\end{tabular}
\end{center}
\caption{Error from truncation in surface harmonic expansion for RBC 1 (50\% weight) sample surface.}
\label{tab:reconstruct2}
\end{table}

\begin{table}[h,t,b]
\begin{center}
\begin{tabular}{lllll} \hline
$N$ & $E_{rms}$ & $E_{Vol}$ & $E_{SA}$ & $E_{Eng}$   \\
\hline
1 & $7.13\times10^{-2}$ & $2.82\times10^{-1}$ & $3.12\times10^{-1}$ & $3.82\times10^{-1}$ \\
2 & $2.06\times10^{-3}$ & $1.35\times10^{-2}$ & $9.60\times10^{-3}$ & $3.96\times10^{-1}$ \\
3 & $2.06\times10^{-3}$ & $1.35\times10^{-2}$ & $9.60\times10^{-3}$ & $3.96\times10^{-1}$ \\
4 & $1.26\times10^{-3}$ & $3.36\times10^{-3}$ & $5.71\times10^{-3}$ & $7.75\times10^{-2}$ \\
5 & $1.26\times10^{-3}$ & $3.36\times10^{-3}$ & $5.71\times10^{-3}$ & $7.75\times10^{-2}$ \\
6 & $3.07\times10^{-4}$ & $1.92\times10^{-3}$ & $5.79\times10^{-3}$ & $4.24\times10^{-1}$ \\
7 & $3.07\times10^{-4}$ & $1.92\times10^{-3}$ & $5.79\times10^{-3}$ & $4.24\times10^{-1}$ \\
8 & $2.15\times10^{-4}$ & $1.21\times10^{-3}$ & $7.14\times10^{-3}$ & $9.97\times10^{-2}$ \\
\hline
\end{tabular}
\end{center}
\caption{Error from truncation in surface harmonic expansion for RBC 2 (90\% weight) sample surface.}
\label{tab:reconstruct3}
\end{table}

\begin{figure}[h!tbs]
    \centering
   \includegraphics[trim=0cm 2cm 0cm 2cm,clip=true,width=8cm]{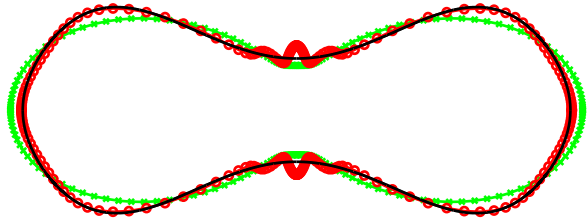} 
\caption{Effect of large $N$ for RBC 1.  The analytical surface is in black, the reconstructed surfaces for $N=4$ and $N=12$ are in green and red, respectively.  Increasing $N$ near 4 improves the pointwise accuracy especially near the edges, but also introduces oscillations near the center with high energy cost.  For $N=12$, the pointwise error in the surface is improved at $1.18\times10^{-4}$, but the error in the energy is 258\%.}  \label{fig:reconstructN} 
\end{figure}

For the red blood cell structures, initially the errors decrease as $N$ increases, but increasing the number of modes beyond a certain threshold actually increases the error in the energy computation.  For RBC 1, the best possible error in the energy is 5.3\%, with $N=4$ or $N=5$ modes.  For RBC 2, the best error is 7.75\% with the same $N$.  We suggest the reason for this is because higher modes contain more bulges than the lower modes, akin to Runge's phenomenon in high order polynomials.  In the reconstruction, the coefficients are chosen to minimize $E_{rms}$.  While transiently increasing $N$ does improve the accuracy of $E_{rms}$, it may introduce local oscillations.  Since the energy is a function of the square mean curvature, these oscillations have a high energy cost.  In Figure \ref{fig:reconstructN}, RBC 1 is reconstructed with $N=4$ and $N=12$, for comparison.  


\section{Examples: reduced volume}

In this section, we provide numerical examples to test the robustness method.  First, observe that the integration of the square local mean curvature $H^2$ is a dimensionless quantity.  When the spontaneous curvature $C_0 = 0$, the mechanical bending energy \eqref{eq:bending_energy} is completely governed by this dimensionless quantity and is therefore scale-invariant.  Thus, for vesicle shapes, where $C_0 = 0$, the minimum energy is completely determined by a single dimensionless quantity called the reduced volume $v$.  If we denote the current vesicle volume and surface area $V$ and $S_A$, respectively, then the reduced volume scales the current volume $V$ by the volume of a sphere with surface area $S_A$.  Since spheres maximize volume for a given surface area, the reduced volume satisfies $v \leq 1$.  The reduced volume is given by the formula
\begin{equation}\label{eq:reduced_v}
v = \DD\frac{V}{4\pi/3R_0^3}
\end{equation}
where $R_0 = \sqrt{S_A/4\pi}$.  In terms of the surface area, the reduced volume is 
\begin{equation}
v = \DD\frac{6\sqrt{\pi}V}{(S_A)^{3/2}}.
\end{equation}

Seifert et. al. has compiled a library of reduced volumes and their corresponding minimum energies for axisymmetric shapes in \cite{Seifert1991} by solving the Euler-Lagrange equations using a parameterization of the vesicle shape with an axis of symmetry.  For verification purposes, we compare our axisymmetric results for various reduced volumes to theirs.  We set the constraint volume $\bar{V}$ to be proportional to $V$ by the (projected) reduced volume $v$.  The volume constraint is in violation and NCG begins to change the shape to relax this configuration.  If we begin with a perfectly spherical vesicle, NCG will simply scale the sphere to a sphere with a smaller volume, and the final reduced volume will be 1.  Therefore, we take a slightly perturbed sphere to be our initial configuration.  After NCG has converged, we calculate the final reduced volume $v$ and the final energy $E$ scaled by the energy of a spherical vesicle $E_0 = 8\pi\mathscr{K}_C$.  

\begin{figure}[h!tbs]
    \centering
   \includegraphics[height=5cm]{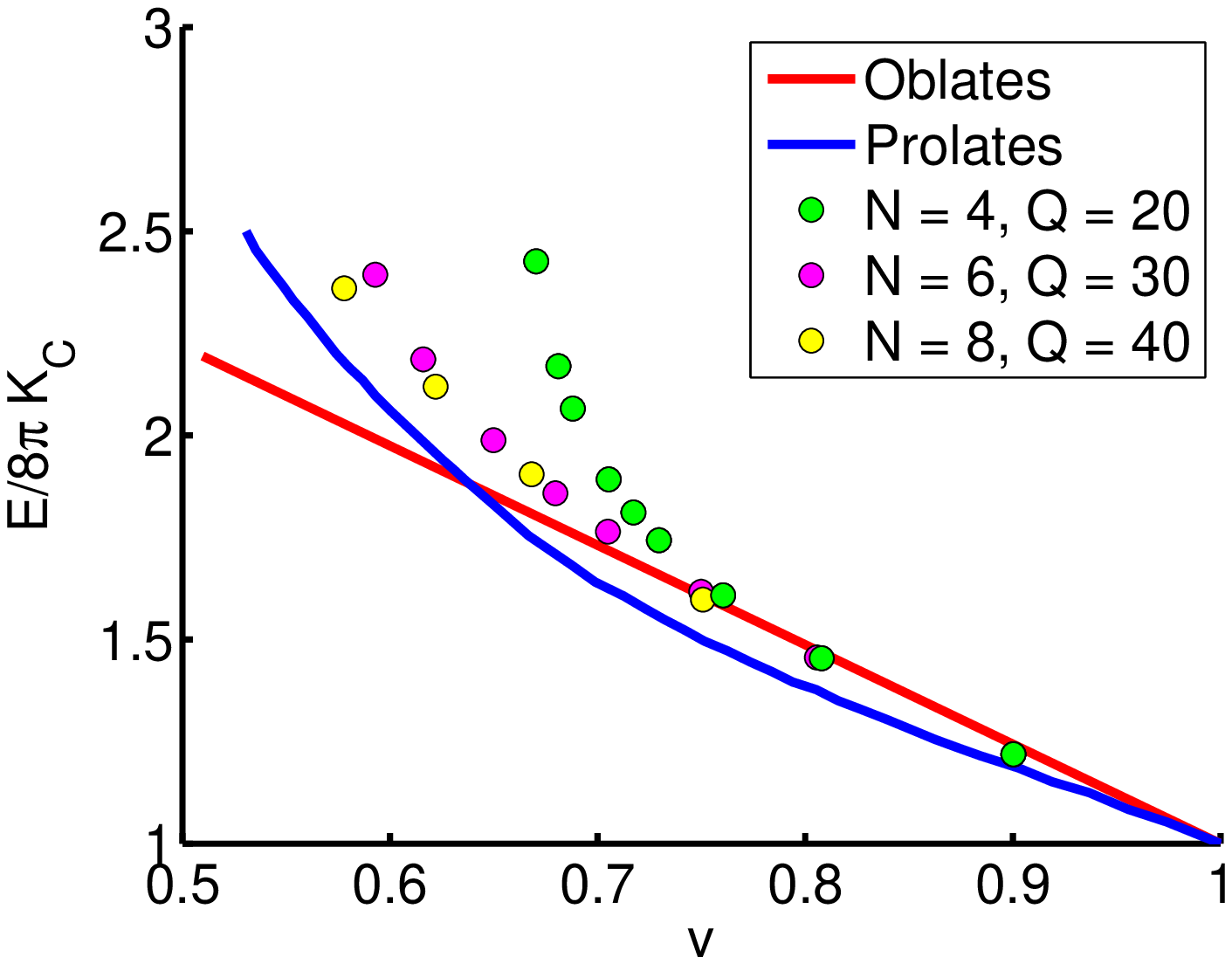} 
   \includegraphics[height=5cm]{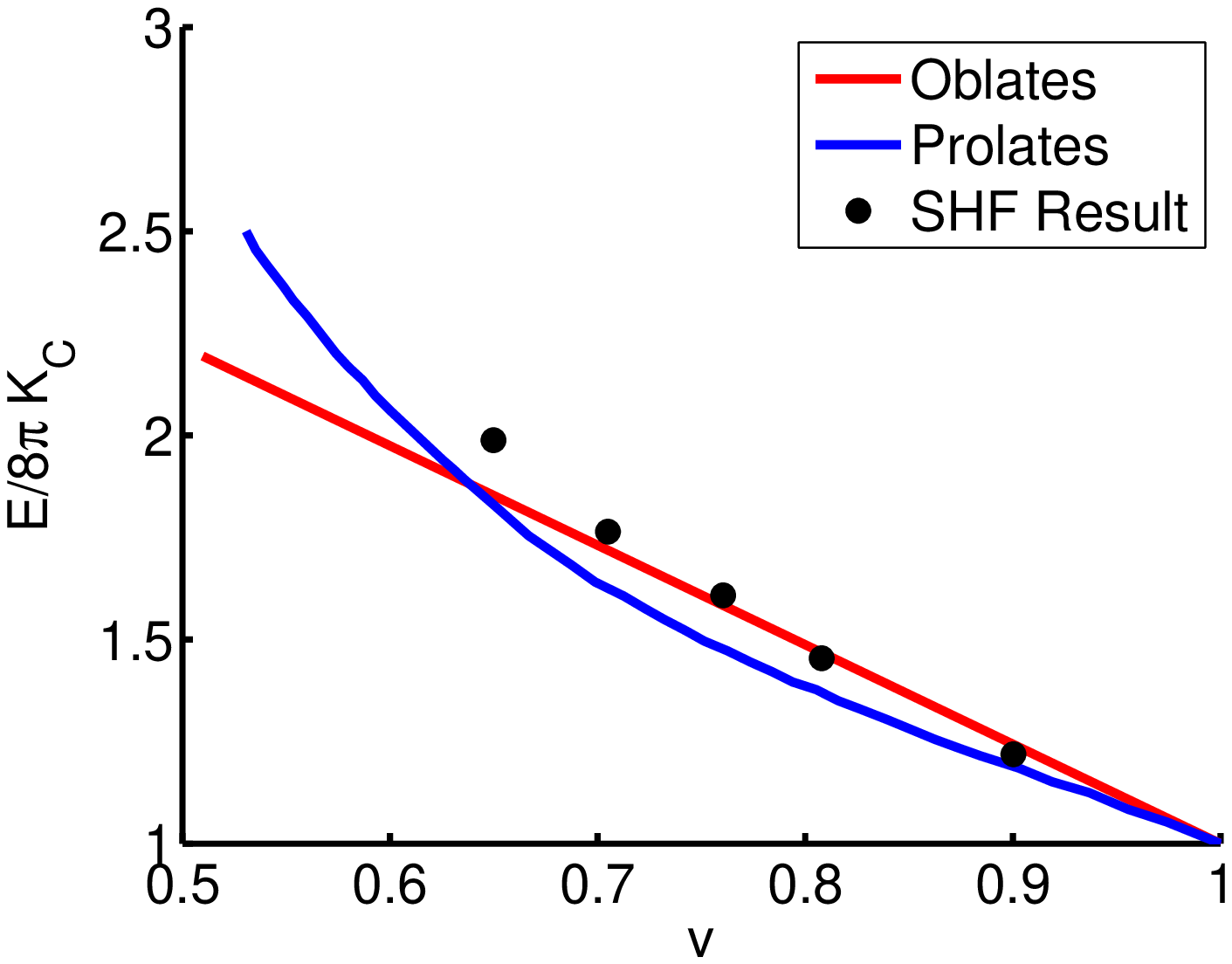}
\caption{Vesicle energy vs. reduced volume for comparison with \cite{Seifert1991}.  Left: All data points used to determine cutoff values for $N$ given a reduced volume $v$.  Right: Kept data according to the cutoff values.  All data points are within 10\% error of the appropriate analytical curve.}  \label{fig:overlay} 
\end{figure}

From the reconstruction examples, $N=4$ is a reasonable truncation for the surface harmonic expansion.  During the iterations of NCG, 20 cubature points are used.  We say that NCG has converged when the norm of the change in gradient or if the change in the modes was less than $10^{-6}$.  When the final configuration is achieved, the total energy is evaluated with 64 cubature points to provide a more accurate computation and to ensure that enough cubature points are used.  For reduced volumes above approximately $v = 0.75$, the numerical energy is within 10\% of the analytical values calculated by Seifert et. al. (see Table \ref{tab:compare}.  However, for reduced volumes less than this, the error exceeds 10\%.  If the number of modes is increased to $N=6$ (since $N=5$ gives the same numerical results as $N=4$ as demonstrated in Tables \ref{tab:reconstruct1}-\ref{tab:reconstruct3}), the relative error is reduced.  However, there is a significant difference between the energy evaluated at 20 cubature points than at 64 cubature points at the final iteration.  This is because the added oscillation from the higher order modes is not absorbed by NCG with only 20 cubature points.  Using 30 cubature points when $N = 6$, the relative error is less than 1\% when compared to 64 points.  For surfaces with reduced volume $0.65 \leq v \leq 0.75$, using $N=6$ and 30 cubature points, the relative error in the final energy is less than 10\%.  For surfaces with reduced volume $0.5 \leq v < 0.65$, we determined that 40 quadrature points are needed to use $N = 8$; however, the error is still above 10\%, and the use of $N = 8$ fared no better than $N=6$.  Our method could not reconstruct surfaces with these reduced volumes well.  These data are plotted in Figure \ref{fig:overlay}, overlayed by the analytic solution from Seifert \cite{Seifert1991}.  

%

\begin{figure}[h!tbs]
    \centering
   \includegraphics[height=2.5cm]{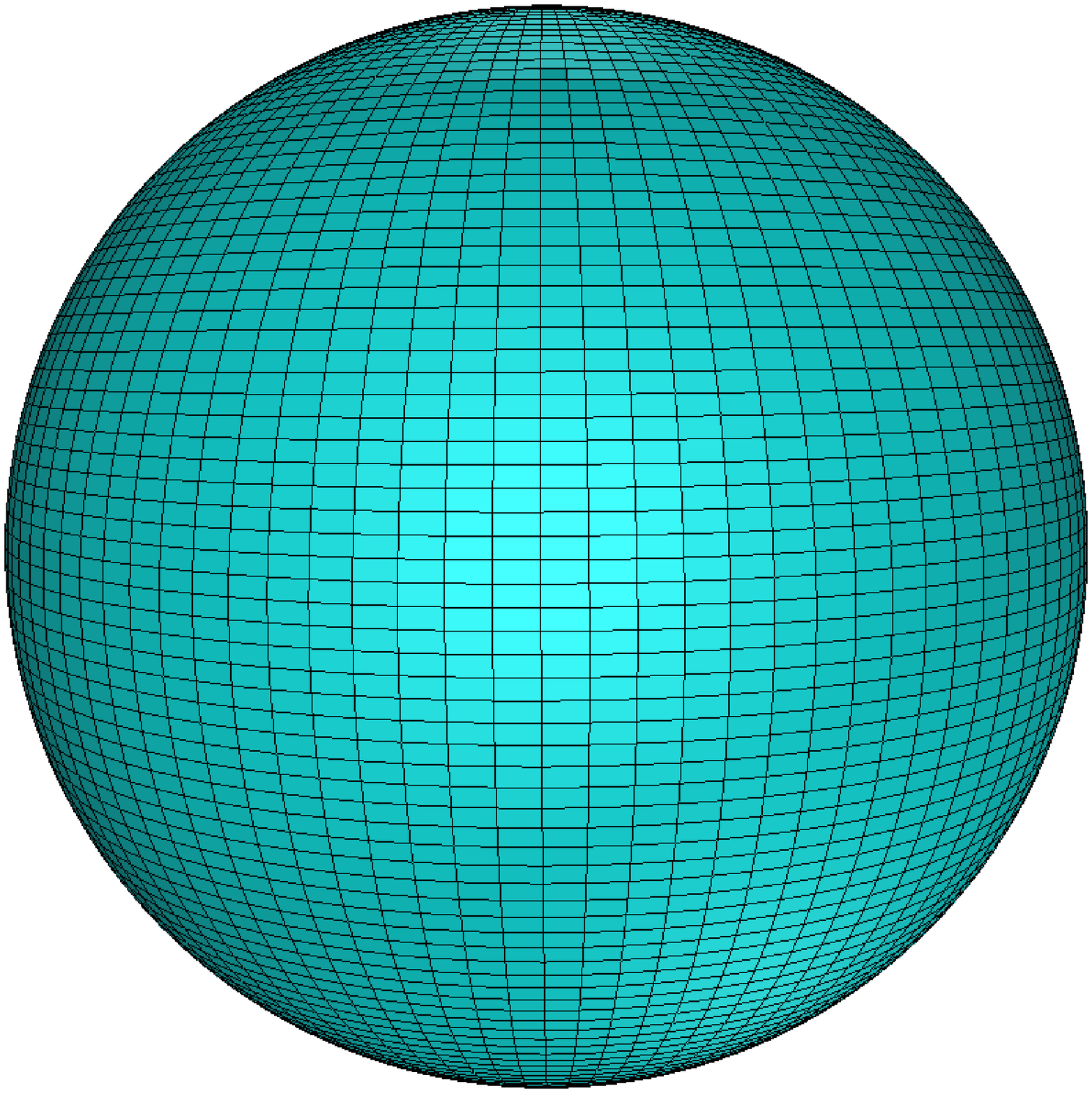}
   \includegraphics[height=2.5cm]{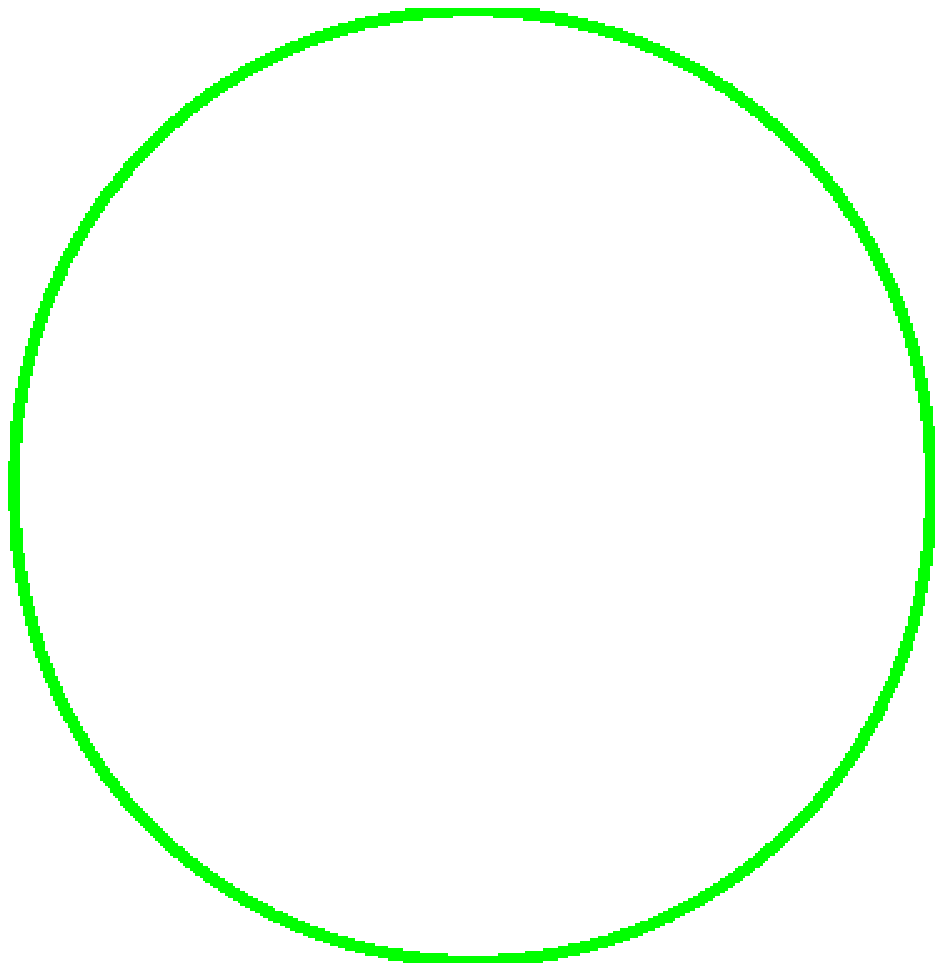} 
   \includegraphics[height=2.5cm]{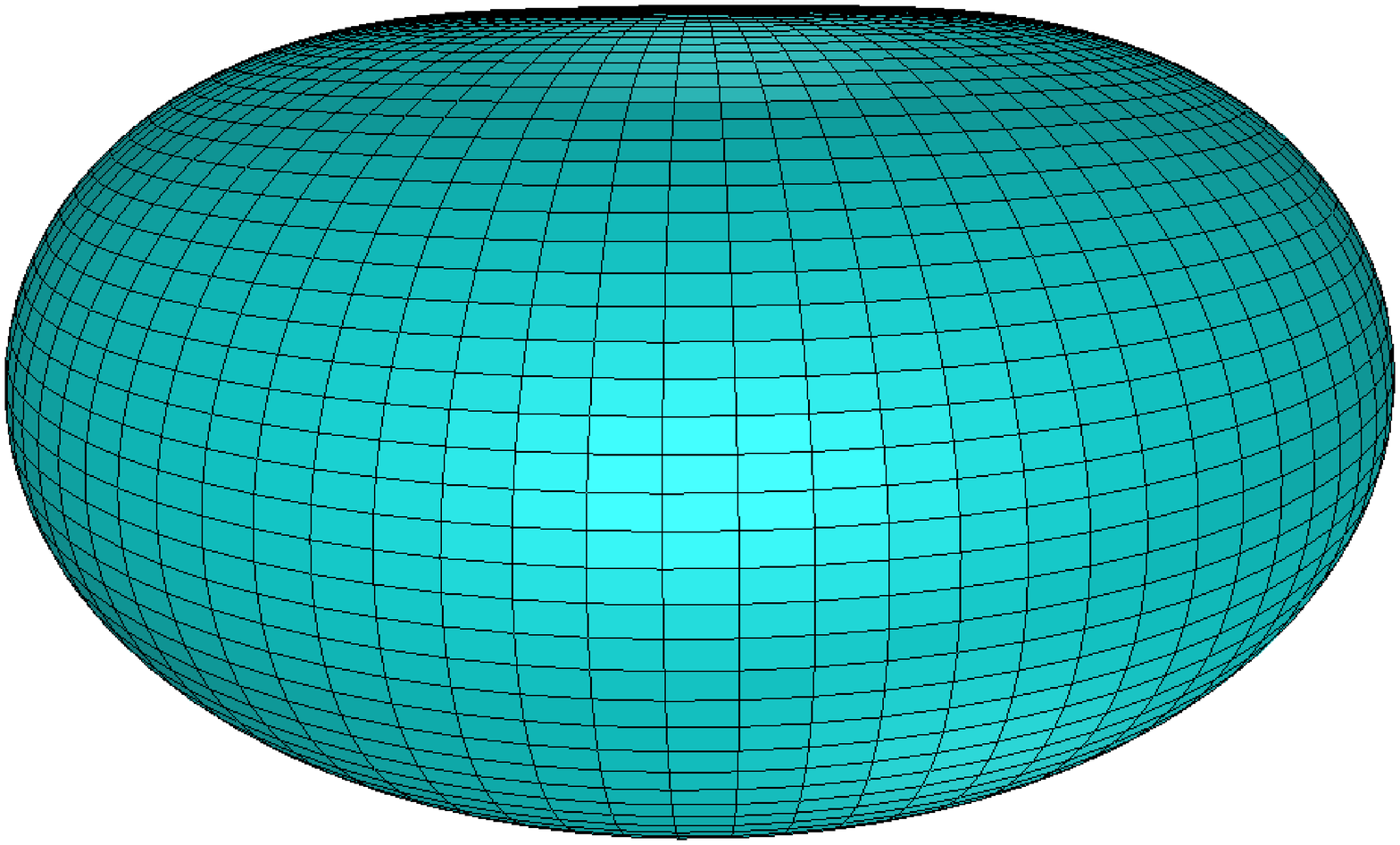}
   \includegraphics[height=2.5cm]{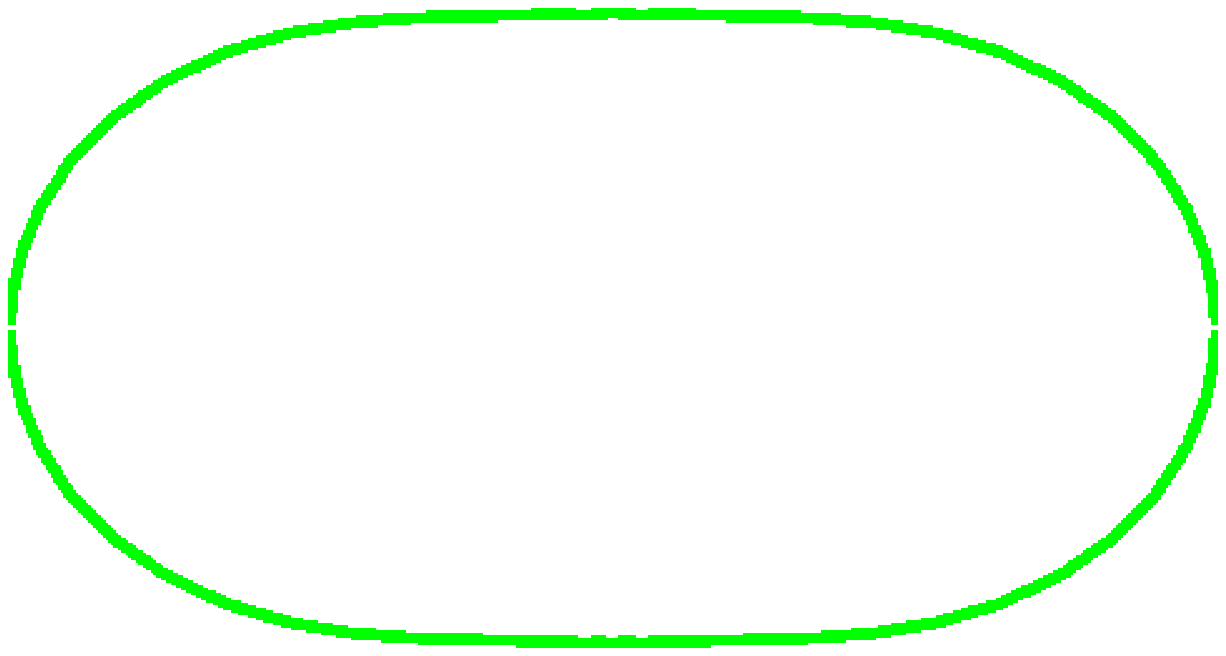} \\
   \includegraphics[height=2.5cm]{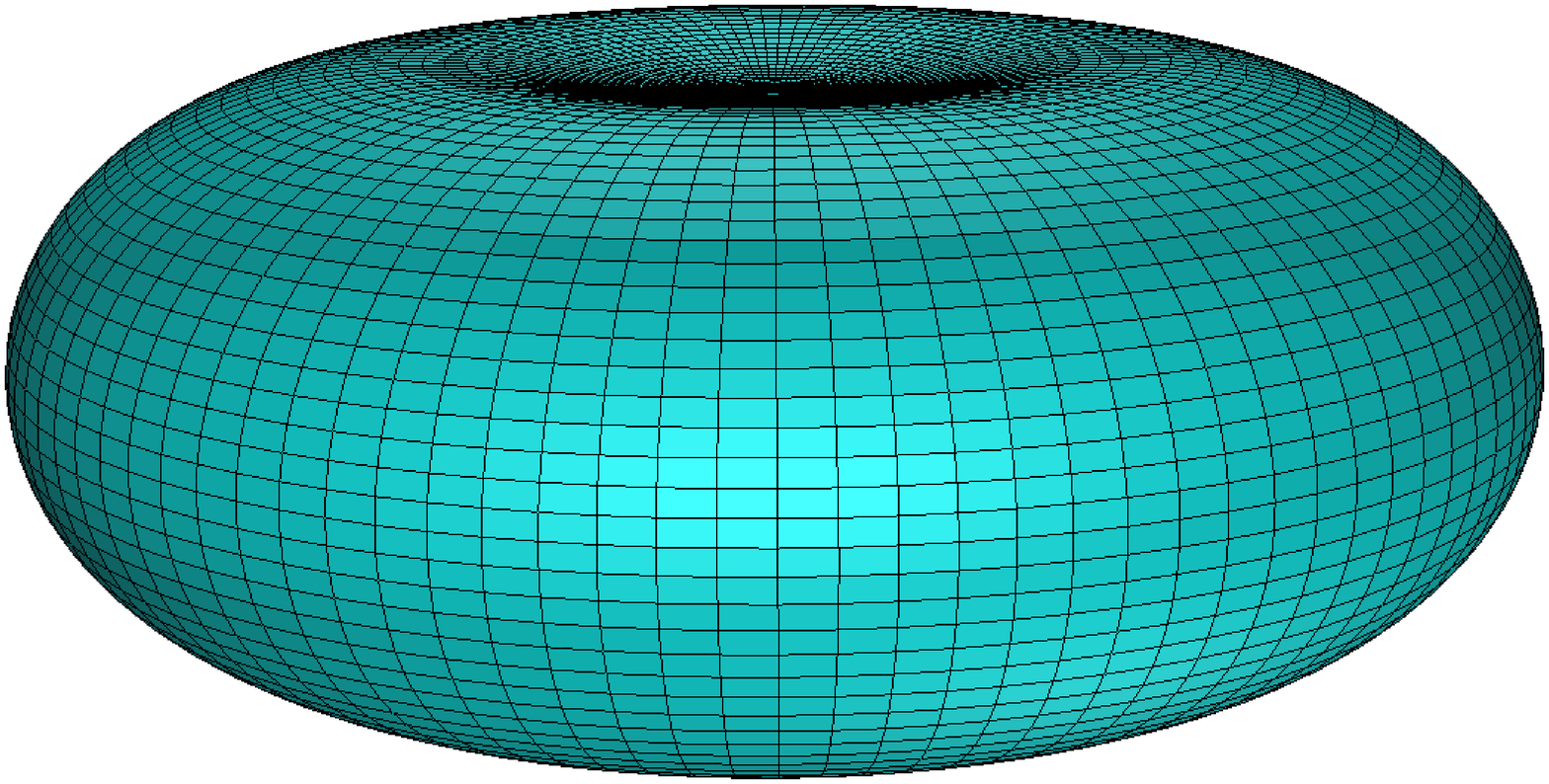}
   \includegraphics[height=2.5cm]{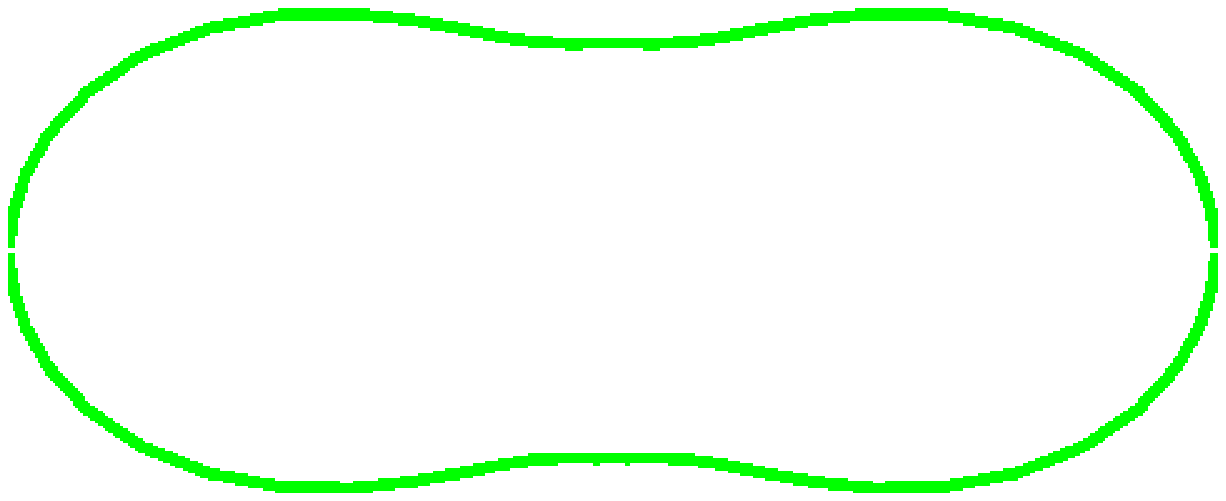}  
   \includegraphics[height=2.5cm]{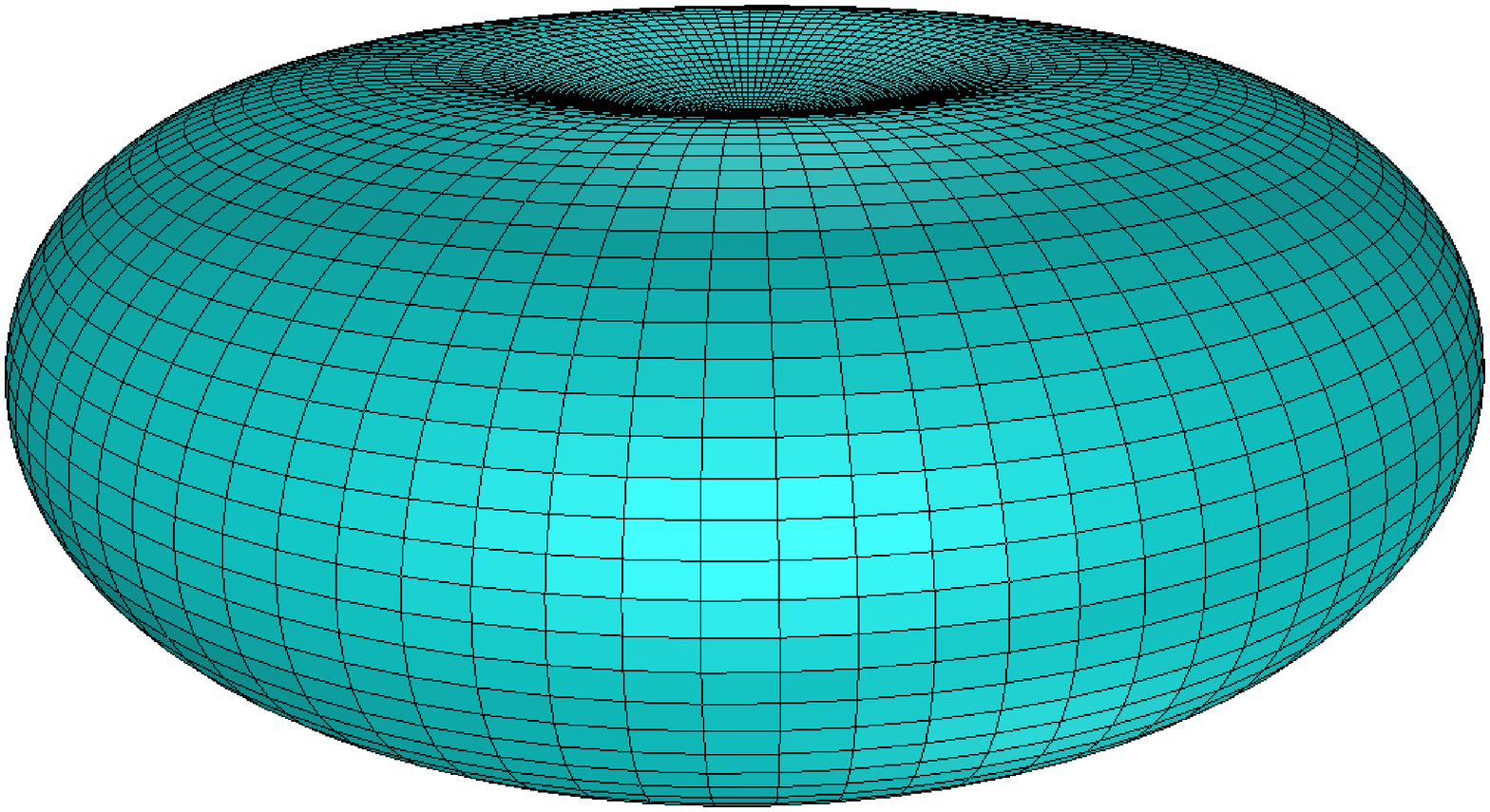}
   \includegraphics[height=2.5cm]{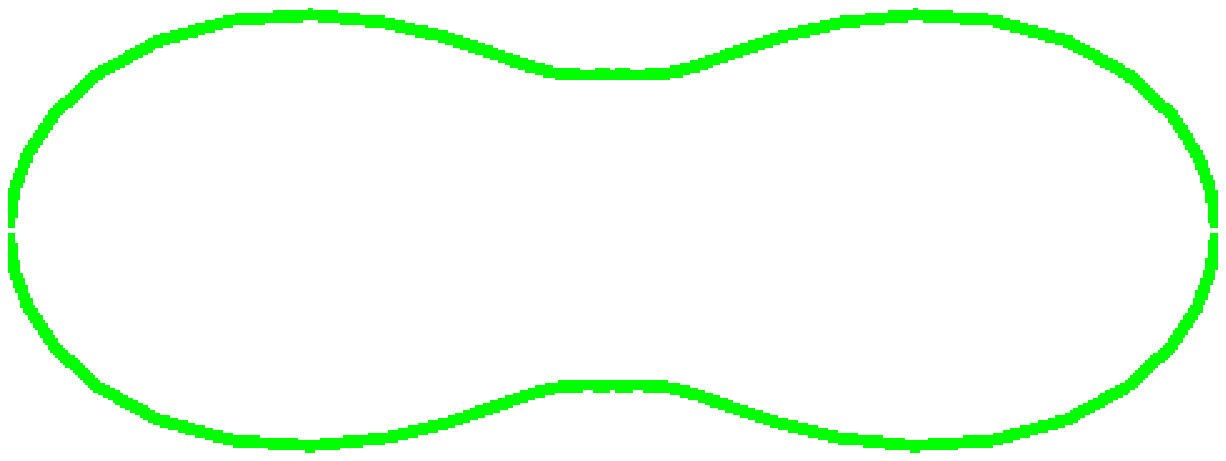} \\ 
   \includegraphics[height=2.5cm]{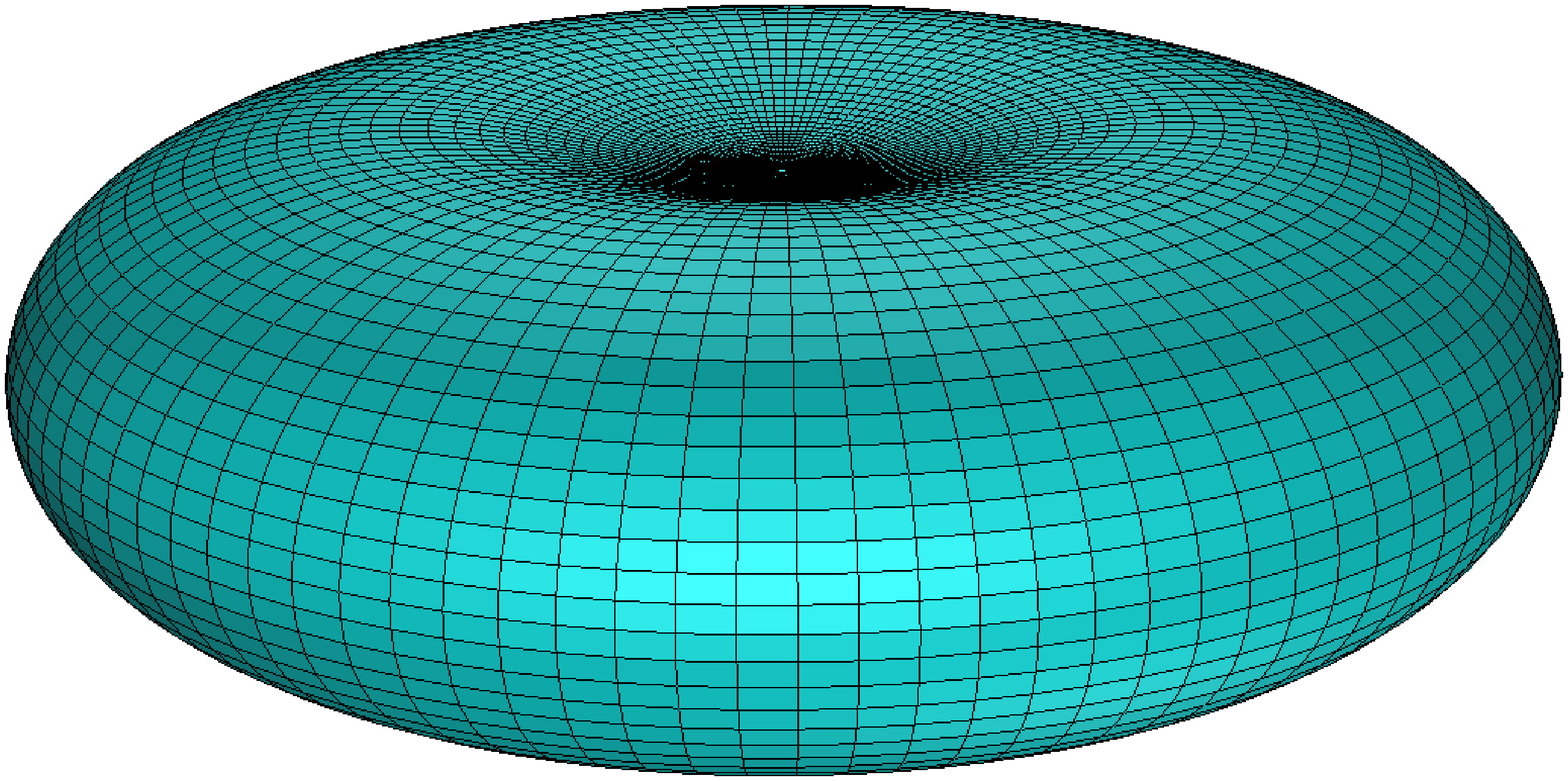}
   \includegraphics[height=2.5cm]{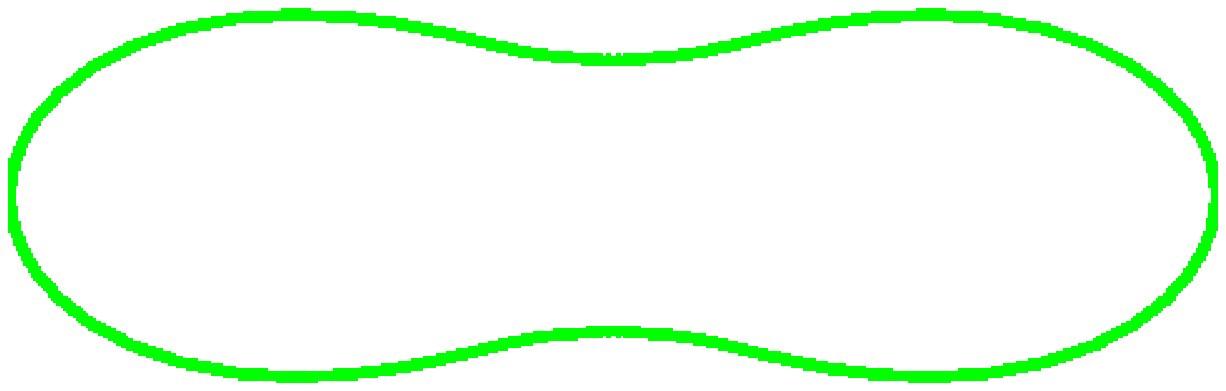} 
   \includegraphics[height=2.5cm]{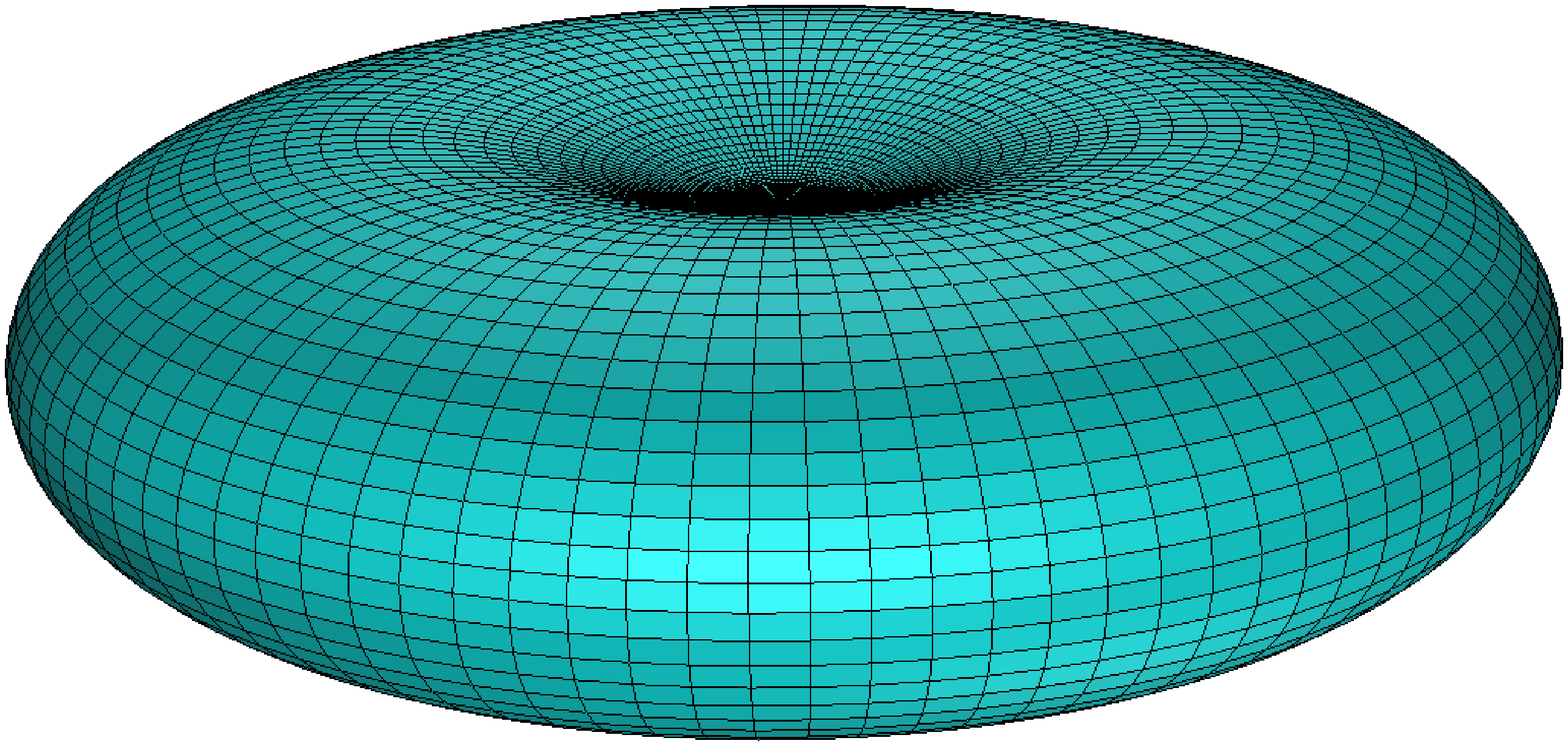} 
   \includegraphics[height=2.5cm]{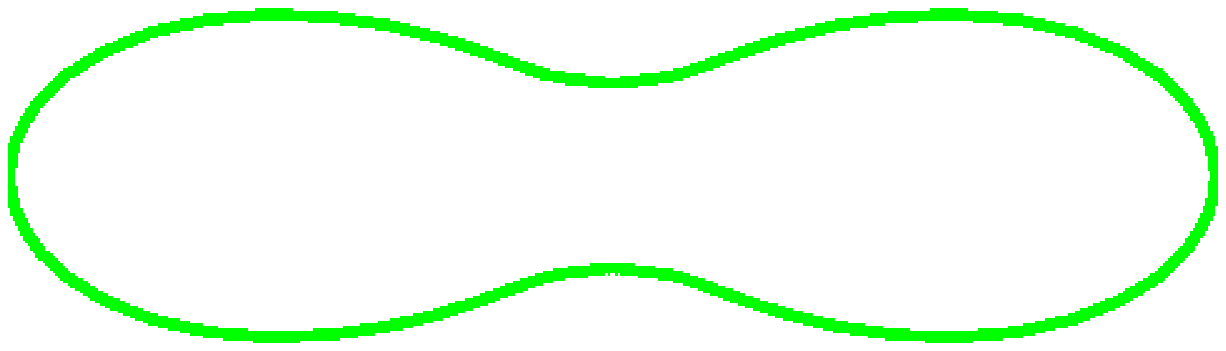} 
\caption{Axisymmetric results with various reduced volumes.  Organized by row, $v = 1.00, 0.90; 0.81, 0.76; 0.70, 0.65$.  The shapes here correspond to the data in Table \ref{tab:compare}.  The full volumes are angled slightly down, while the profiles are not angled.}  \label{fig:shapes} 
\end{figure}

In summary, for surfaces with $0.75 \leq v \leq 1$, we used $N = 4$ and 20 cubature points, for surfaces with $0.65 \leq v < 0.75$, we used $N=6$ and 30 cubature points.  The results using this cutoff are overlayed by Seifert's data in Figure \ref{fig:overlay}.  Finally, the surfaces corresponding to the data in Table \ref{tab:compare} are shown in Figure \ref{fig:shapes}.  

\begin{table}[ht]
\begin{center}
\begin{tabular}{lllllll} \hline 
$v$                         & 1.0 & 0.90  & 0.81  & 0.76  & 0.70  & 0.65  \\
\hline
$E/E_0$, \cite{Seifert1991} & 1.0 & 1.19  & 1.37  & 1.58  & 1.72  & 1.85  \\

$E/E_0$, SHF                & 1.0 & 1.22  & 1.45  & 1.61  & 1.76  & 1.99  \\

Rel. error                  & 0\% & 2.5\% & 5.8\% & 1.9\% & 2.3\% & 7.6\% \\
\hline
\end{tabular}
\end{center}
\caption{Vesicle energy vs. reduced volume for comparison with \cite{Feng2006}.}\vspace{0.5em}
\label{tab:compare}
\end{table}

We note that the results of the numerical procedure may be only local minima and therefore only locally stable.  With enough perturbation through some external force, another configuration with a lower energy may be achieved.  In the range of $0.64 \leq v \leq 1$, oblate shapes are local energy minimizers, but prolates are global minimizers for axisymmetric shapes.  However, it may be possible to obtain a non-axisymmetric shape with lower energy than a prolate.  

\section{Conclusion}
We have presented a fast algorithm for computing axisymmetric and non-axisymmetric spherical vesicles that correspond to 
minimized Canham-Helfrich-Evans curvature energy subject to surface area and volume constraints. Our method is based  
on the real-valued surface harmonic rather than the complex-valued spherical harmonic expansion of the surface configurations.
Computational simulations showed that vesicles of various reduced volumes can be approximated using up to 49 surface
harmonic functions, and the approximation error measured in curvature energy can be well maintained within 8\%, mostly
below 5\% indeed. We use a nonlinear conjugate gradient method rather than the Newton's method for the numerical minimization. 
This makes it possible for us to avoid the computation of the Hessian matrices and the solution of linear systems to
further improve the efficiency of the calculation. Our method entails advantageous over the spherical harmonic approximations
of the individual coordinates of the vesicle configurations since it excludes the complex parts unnecessary for computing real-valued 
surfaces \cite{KhairyK2011a}. We will explore the implementation of this fast algorithm in the external
force fields such as fluid flow or electrostatic potential field, for which re-orientation of the vesicle
configurations might be necessary because of the vesicle rotation caused by the non-vanishing torque applied by
the external force. A fast algorithm for the vesicle deformation induced by the electrostatic force is currently
under development and will be reported elsewhere.

\appendix
\section{Nonlinear Conjugate Gradient Method and Linear Search} \label{append}
The pseudocode of the nonlinear conjugate method and the related linear search method is given below
to ease the implementation. The functions $I(\cdot)$ and $\delta_{\Gamma} I(\cdot)$ evaluate the curvature energy and its
gradient, c.f. Eqs. (\ref{eq:I}) and (\ref{eq:dI}), respectively.
\begin{framed}
\noindent \textbf{Algorithm 3.1:} \textit{Nonlinear Conjugate Gradient} \\
\indent Given initial SHF modes $\vec{a}_0$, tolerances $\epsilon_g, \epsilon_{a}$, $M$ \\
\indent Compute $E_0 \leftarrow I(\vec{a}_0)$, $\vec{g}_0 \leftarrow \delta_\Gamma I(\vec{a}_0)$ \\
\indent Direction $\vec{d}_0 = -\vec{g}_0$\\
\indent $k \leftarrow 0$ \\
\indent \textbf{for} $k = 1: M$ \\
  \indent \qquad Step size $\alpha_k \leftarrow$ \textit{LineSearch}$(\vec{a}_k, \vec{d}_k)$ \\
  \indent \qquad Step in direction $\vec{a}_{k+1} \leftarrow \vec{a}_k + \alpha_k \vec{d}_k$ \\
  \indent \qquad Update energy and gradient $E_{k+1} \leftarrow I(\vec{a}_{k+1})$, $\vec{g}_{k+1} \leftarrow \delta_\Gamma I(\vec{a}_{k+1})$ \\
  \indent \qquad \textbf{if} $|| \vec{g}_{k+1} - \vec{g}_k|| < \epsilon_g $ ~ \textbf{break} \\
  \indent \qquad Compute $\beta_k \leftarrow (\vec{g}_{k+1}^T(\vec{g}_{k+1}-\vec{g}_k)) / ((\vec{g}_{k+1} - \vec{g}_k)^T\vec{d}_k)$ \\
  \indent \qquad Update direction $\vec{d}_{k+1} \leftarrow - \vec{g}_{k+1} + \beta_k \vec{d}_k$ \\
  \indent \qquad \textbf{if} $|| \vec{a}_{k+1} - \vec{a}_k || < \epsilon_a$ ~ \textbf{break} \\
  \indent \qquad $k \leftarrow k+1$\\
  \indent \textbf{end} \\
\indent \textbf{return} $\vec{a}_k$ 
\end{framed}

\begin{framed}
\noindent \textbf{Algorithm 3.2:} \textit{LineSearch}$(\vec{a}, \vec{d})$ \\
\indent Define $\alpha_m \leftarrow 0, \alpha_M \leftarrow 1, M, \epsilon$ \\
\indent Compute $E_m \leftarrow I(\vec{a} + \alpha_m \cdot \vec{d})$,  
$E_M \leftarrow I(\vec{a} + \alpha_M \cdot \vec{d})$ \\
\indent \textbf{if} $E_m < E_M$ \\
  \indent \qquad $\alpha_l \leftarrow \alpha_m$; $\alpha_u \leftarrow \alpha_M$ \\
  \indent \qquad $E_l \leftarrow E_m$; $E_u \leftarrow E_M$ \\
\indent \textbf{else} \\
  \indent \qquad $\alpha_l \leftarrow \alpha_M$;  $\alpha_u \leftarrow \alpha_m$ \\
  \indent \qquad $E_l \leftarrow E_M$; $E_u \leftarrow E_m$ \\
\indent \textbf{end if} \\
\indent \textbf{for} $i = 1:M$ \\
  \indent \qquad $\alpha_t \leftarrow (\alpha_l + \alpha_u)/2$ \\
  \indent \qquad $E_t \leftarrow I(\vec{a}+ \alpha_t \cdot \vec{d})$ \\ 
  \indent \qquad \textbf{if} $E_t > E_l$ \textbf{then} \\
    \indent \qquad \qquad  $\alpha_u \leftarrow \alpha_t$ \\
  \indent \qquad \textbf{else} \\
    \indent \qquad \qquad $\vec{g} \leftarrow \delta_\Gamma I(\vec{a} + \alpha_t \cdot \vec{d})$ \\
    \indent \qquad \qquad $(D\phi) \leftarrow \vec{g}^T \, \vec{d}$ \\
    \indent \qquad \qquad \textbf{if} $(D\phi) \cdot (\alpha_l - \alpha_t) >0$ \textbf{then} \\
      \indent \qquad \qquad \qquad $\alpha_l \leftarrow \alpha_t$ \\
    \indent \qquad \qquad \textbf{else} \\
      \indent \qquad \qquad \qquad $\alpha_u \leftarrow \alpha_l$; $\alpha_l \leftarrow \alpha_t$ \\
    \indent \qquad \qquad \textbf{end if} \\
    \indent \qquad \qquad $E_l \leftarrow E_t$ \\
  \indent \qquad \textbf{end if} \\
  \indent \qquad \textbf{if} $|\alpha_u - \alpha_l| < \epsilon$ ~ \textbf{break} \\
\indent \textbf{return} $\alpha_t$
\end{framed}

\newpage

\end{document}